\documentclass[oneside, 10pt]{amsart}

\usepackage[a4paper, margin=1.4in, bottom=2in, footskip=0.8in]{geometry} 
\usepackage{setspace} 
\usepackage{appendix}

\usepackage{setspace}

\usepackage{pgfplots}
\usepackage{float} 
\pgfplotsset{compat=newest}
\usetikzlibrary{intersections}
\usepgfplotslibrary{fillbetween}

\usepackage[style=numeric, sorting=none,backend=biber]{biblatex}
\usepackage{hyperref}
\hypersetup{
    colorlinks = true, 
    urlcolor = blue, 
    linkcolor = blue, 
    citecolor = blue 
}

\usepackage{amsmath}
\usepackage{amsthm}
\usepackage{textcomp}

\newtheorem{theorem}{THEOREM}[section]
\newtheorem{lemma}[theorem]{LEMMA}
\newtheorem{proposition}[theorem]{PROPOSITION}
\newtheorem{corollary}[theorem]{COROLLARY}
\theoremstyle{definition}
\newtheorem{definition}[theorem]{DEFINITION} 
\theoremstyle{remark}

\newtheorem*{remark}{Remark}

\newcommand{\thesistitle}[0]{The global existence of small-amplitude solutions to nonlinear Klein-Gordon equations: A study based on S. Klainerman's approach}
\newcommand{\authorname}[0]{Alessandro Massaad}

\newcommand{\supervisor}[0]{Annalaura Stingo}
\newcommand{\supervisorinstitution}[0]{Centre de Mathématiques Laurent Schwartz }

\usepackage[T1]{fontenc}
\usepackage{amstext}
\usepackage{amsmath}
\usepackage{amssymb}
\usepackage{url}
\usepackage{graphicx}
\usepackage{wrapfig}
\usepackage{enumerate}
\usepackage{paralist}
\usepackage{xspace}
\usepackage{color}
\usepackage{times}

\usepackage{fancyhdr} 

\usepackage{titling}

\usepackage{amsthm}

\newcommand{\lcomb}[1]{\sum^\prime \limits_{#1}{\vphantom{\sum}}}

\begin{document}

\onehalfspacing


\thispagestyle{empty}


\hspace{0pt}
\vfill

\begin{center}

\includegraphics[width=0.3\textwidth]{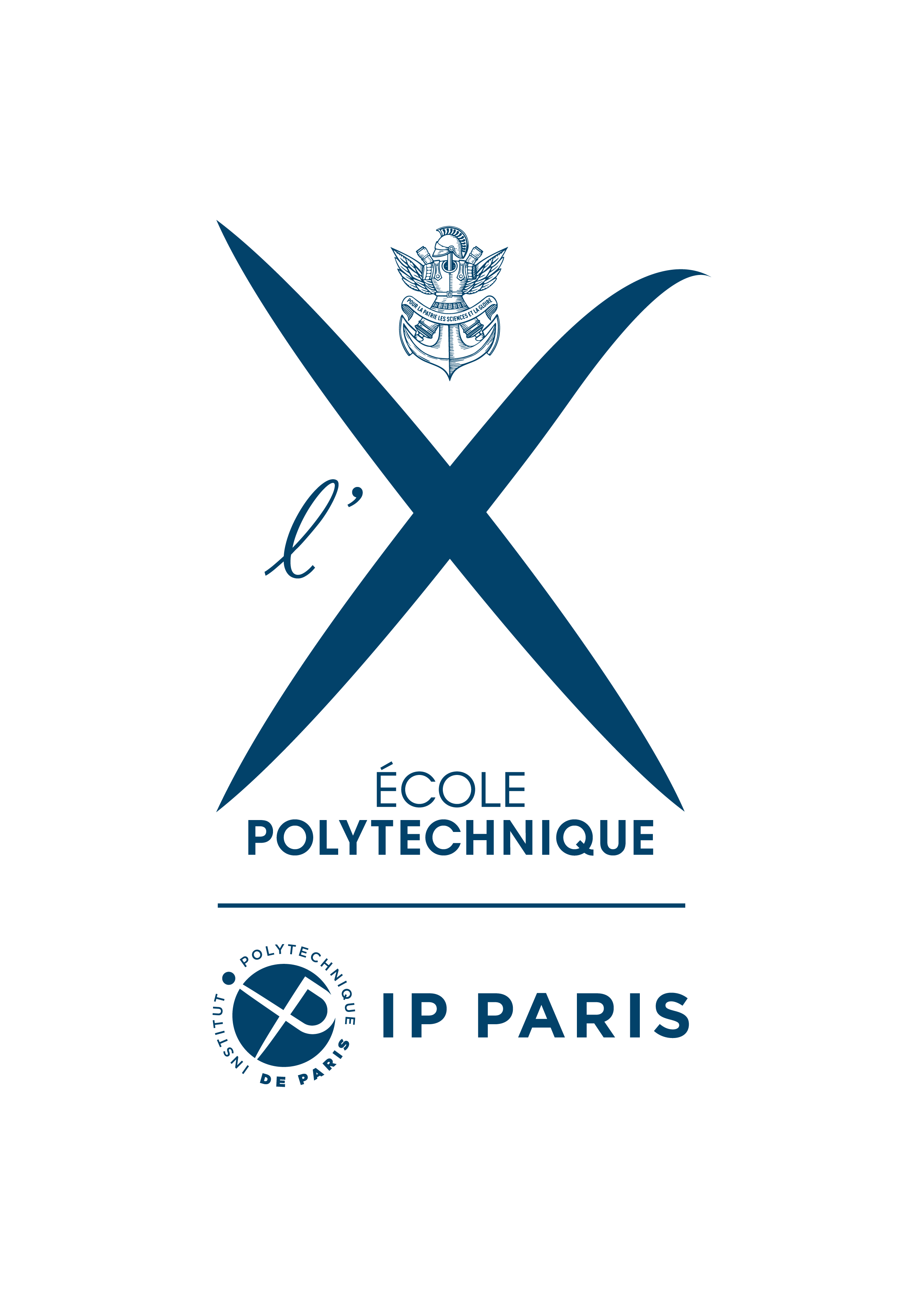}

\vspace*{0em}
{\large 

\vspace*{0em}
\textit{BACHELOR THESIS IN MATHEMATICS}

\vspace*{3em}
{\Huge \textbf{\thesistitle}}
\vspace*{6em}

\textit{Author:}

\vspace*{1em}
\authorname{}, \'Ecole Polytechnique

\vspace*{2em}
{\textit{Advisor:}}

\vspace*{1em}
\supervisor{}, \supervisorinstitution{}
}

\vspace*{10em}
\textit{March 2024}
\end{center}

\newpage

\pagestyle{fancy}
\lhead{\normalsize \textit{ALESSANDRO MASSAAD}}
\rhead{\normalsize \textit{BACHELOR THESIS}}

\newpage
\pagenumbering{arabic}

\section*{Abstract}
In this thesis we explore S. Klainerman's proof on the global existence of small amplitude solutions to nonlinear Klein-Gordon equations in four space-time dimensions, as established in his paper from 1985 \cite{klainerman1985}. We consider initial data with small amplitude and compact support and aim prove the global existence and uniform decay of smooth solutions. We establish that solutions exist globally if the initial data satisfy a suitable smallness condition. Key analytical tools include generalized Sobolev norms and uniform decay estimates for the associated linear problem. The solutions exhibit a decay rate of $t^{-5/4}$, uniform in time and space. This result is achieved by combining the energy method, perturbed Klein-Gordon techniques, and Sobolev inequalities.

\tableofcontents

\section{Introduction}\label{section:introduction}

\subsection{The problem at hand }
\mbox{}\\
Consider the nonlinear Klein-Gordon equation
\begin{equation}
\Box u+u=F\left(u, u^{\prime}, u^{\prime \prime}\right)\tag{N.K.G.}\label{eq:NKG}
\end{equation}

where $\square=\partial_{t}^{2}-\partial_{1}^{2}-\partial_{2}^{2}-\partial_{3}^{2}$ is the D'Alembertian of 4-dimensional space-time $\mathbb{R}_+\times\mathbb{R}^3$, $u : \mathbb{R}_+\times\mathbb{R}^3\rightarrow\mathbb{R}$ is a scalar function and $F$ a smooth function of $u$, $u'$ (representing all of the first partial derivatives of $u$) and $u''$ (representing all of the second partial derivatives of $u$). Assume that $F$ vanishes, together with its first partial derivatives, at $\left(u, u^{\prime}, u^{\prime \prime}\right)= (0, 0, 0)$. In other words, $F$ is at least quadratic in $u, u^{\prime}, u^{\prime \prime}$ in a neighborhood of the trivial solution $u \equiv 0$. We subject $u$ to the initial value problem
\begin{equation}
\left\{
\begin{aligned}
u(0, \;\cdot\;)&=\varepsilon f, \\
\partial_t u(0, \;\cdot\;)&=\varepsilon g
\end{aligned}
\right.\tag{I.V.P.}\label{eq:IVP}
\end{equation}
with $f, g \in C_{0}^{\infty}\left(\mathbb{R}^{3}\right)$ and $\varepsilon>0$, a small parameter. In this thesis, we aim to prove the following result, using S. Klainerman's approach in \cite{klainerman1985}.
\begin{theorem}\label{thm:1}
There exists $\varepsilon_{0}>0$, depending on a finite number of derivatives of $f, g, F$ such that, for any $\varepsilon\in\,]0, \varepsilon_{0}[\,$, \eqref{eq:NKG} has a unique solution, $u \in$ $C^{\infty}\left(\mathbb{R}_+ \times \mathbb{R}^{3}\right)$ which satisfies \eqref{eq:IVP}. Moreover, $u$ decays, uniformly in $t \geqq 0, x \in \mathbb{R}^{3}$, as $t^{-5 / 4}$ for large $t$.
\end{theorem}

The main analytic tools used in our proof are the uniform decay estimates for solutions to the linear inhomogeneous Klein-Gordon equations, $\Box_{1} u=g$, which we develop in Section ~\eqref{section:decay-estimates}. In Section ~\eqref{section:hormander}, we prove the energy method for the perturbed Klein-Gordon equation, a key result in the proof of the main theorem. 

\subsection{The $\Gamma$ Operators}
\mbox{}\\
Our proof of the theorem relies on the Lorentz invariance properties of the linear part $\Box_{1}:=\Box+1$. It is easy to show that the linear homogeneous Klein-Gordon equation $$\Box_1 u = 0$$ is invariant with  respect to translations in time, translations in space, Euclidean rotations (rotations in the spacial domain), and hyperbolic rotations (rotations in the space-time domain). To each of these symmetries we associate a differential operator. The translations are associated to the usual partial derivatives $\{\partial_t, \partial_1, \partial_2, \partial_3\}$. For the rotations, we consider the operators
\begin{align*}
\Omega_{a b}=x_{a} \partial_{b}-x_{b} \partial_{a}, \quad \quad \quad 0 \leq a<b \leq 3,
\end{align*}
where
\begin{align*}
&x_{0}=t, \\
&\partial_{0}=-\partial_{t}, \\ &\partial_{i}=\frac{\partial}{\partial x_{i}}, \quad i=1,2,3
\end{align*}
Thus, 
\begin{align*}
& \Omega_{i j}=x_{i} \partial_{j}-x_{j} \partial_{i}, \\
& \Omega_{0 i}=t \partial_{i}+x_{i} \partial_{t} .
\end{align*}
We associate the Euclidean rotations to the operators $\{\Omega_{i j}\mid 1 \leq i<j \leq 3 \}$, and the hyperbolic rotations to the operators $\{\Omega_{0 i}\mid 1 \leq i \leq 3 \}$. The family $\Gamma := \{\partial_a\mid 0\leq a \leq 3\} \cup \{\Omega_{a b}\mid 0 \leq a<b \leq 3 \}$ generates a Lie algebra, i.e., its $\mathbb{R}$-linear span is a Lie algebra. We define the Lie bracket (the commutator) of any two operators $A,B$ in this Lie algebra by:
\begin{align*}
    [A,B] := AB-BA
\end{align*}
The Lorentz invariance of the Klein-Gordon operator $\Box_{1}=\Box+1$ is then equivalent to
\begin{align*}
\left[\Omega_{a b}, \square_{1}\right]=0, \quad \quad 0 \leq a<b \leq 3
\end{align*}
(see lemma ~\eqref{lemma:lorentz-invatiance})
\subsection{Generalized Sobolev Norms} 
\mbox{}\\
Given a \hyperref[def:multi-index]{10-index $\alpha$}, and for any fixed ordering of the operators $\Gamma$, we denote by $\Gamma^{\alpha}$ the product 
\begin{align*}
    \Gamma^\alpha = \prod_{i=1}^{10}\Gamma_{i}^{\alpha_{i}},
\end{align*}
where $\alpha_i$ is the $i$-th coordinate of $\alpha$ and $\Gamma_i$ is the $i$-th operator in the ordered $\Gamma$ family. 
We can hence define some energy norms for functions $u:\mathbb{R}_+\times\mathbb{R}^3\rightarrow\mathbb{R}$:
\begin{align}\label{def:energynorm}
&\mathcal{E}(u(t, \;\cdot\;)) := \left(\|u(t, \;\cdot\;)\|_{L^2}^2 + \sum_{i=0}^3\|\partial_i u(t, \;\cdot\;)\|_{L^2}^2 \right)^\frac{1}{2}, \\
&\mathrm{E}(u(t, \;\cdot\;)) := |u(t, \;\cdot\;)|_{L^\infty} + \sum_{i=0}^3 |\partial_i u(t, \;\cdot\;)|_{L^\infty}
\end{align}
We can also define the generalized \hyperref[def:sobolev-norms]{Sobolev norms}:
\begin{align}\label{def:generalized-sobolev-norms}
&\|u(t, \;\cdot\;)\|_{\Gamma, N} := \sum_{|\alpha| \leqq N} \mathcal{E}(\Gamma^{\alpha} u(t, \;\cdot\;)) \\
& \quad \quad \quad = \sum_{|\alpha| \leqq N} \left(\|\Gamma^{\alpha} u(t, \;\cdot\;)\|_{L^2}^2 + \sum_{i=0}^3\|\partial_i \Gamma^{\alpha} u(t, \;\cdot\;)\|_{L^2}^2 \right) ^{1 / 2}, \notag \\
&|u(t, \;\cdot\;)|_{\Gamma, N} := \sum_{|\alpha| \leqq N} \mathrm{E}(\Gamma^{\alpha} u(t, \;\cdot\;)) \\
& \quad \quad \quad = \sum_{|\alpha| \leqq N} \left( |\Gamma^{\alpha} u(t, \;\cdot\;)|_{L^\infty} + \sum_{i=0}^3 |\partial_i \Gamma^{\alpha} u(t, \;\cdot\;)|_{L^\infty} \right) \notag
\end{align}
where $\|\cdot\|_{L^2}$ and $|\cdot|_{L^\infty}$ are respectively the $L^2$-norm and the $L^\infty$-norm over the spacial domain $\mathbb{R}^3$. Since the commutator of any two $\Gamma$ operators is a $\mathbb{R}$-linear combination of the $\Gamma$ operators (see lemma ~\eqref{lemma:commutator-of-Gammas}), any two different orderings of the $\Gamma$ 's will produce equivalent norms defined above. We shall also use the following space-time norms for functions $u, g:\mathbb{R}_+\times\mathbb{R}^3\rightarrow\mathbb{R}$:
\begin{align}\label{def:normkN}
&\|u\|_{N}=\sup _{t \geqq 0}\|u(t, \;\cdot\;)\|_{\Gamma, N}, \\
&|u|_{k, N}=\sup _{t \geqq 0}(1+t)^{k}|u(t, \;\cdot\;)|_{\Gamma, N}, \\
&E_{k, N}(g)=\sup _{t \geqq 0}(1+t)^{k}\|g(t, \;\cdot\;)\|_{\Gamma, N},
\end{align}
The smallness condition implicit in the statement of our theorem can be made more precise. Thus, we replace ~\eqref{eq:IVP} by
\begin{align}
\begin{cases}
u(0, \;\cdot\;)&= \quad u_0 \\
\partial_t u(0, \;\cdot\;)&= \quad u_1
\end{cases}\tag{I.V.P.'}\label{eq:IVP'}
\end{align}
with $u_{0}, u_{1} \in C_{0}^{\infty}\left(\mathbb{R}^{3}\right)$ and define $\Delta_{N}=\Delta_{N}\left(u_{0}, u_{1}\right)$ as the restriction of $\|u(t, \;\cdot\;)\|_{\Gamma, N}$ at $t=0$, i.e.,
\begin{equation}\label{def:Delta}
\Delta_{N}=\left.\left(\|u(t, \;\cdot\;)\|_{\Gamma, N}\right)\right|_{t=0}
\end{equation}
$\Delta_{N}$ depends of course, only on $u_{0}, u_{1}$. The conclusion of the theorem holds if $\Delta_{N_{0}} \leqq \varepsilon$ for any $0<\varepsilon \leqq \varepsilon_{0}$, with $N_{0}$ a fixed positive integer to be specified during the proof.

\section{Decay Estimates for the Linear Inhomogenous Problem}\label{section:decay-estimates}

Consider the linear Klein-Gordon equation
\begin{equation}\label{eq:LKG}
\Box_{1} u=g\tag{L.K.G.}
\end{equation}
subject to the initial conditions \eqref{eq:IVP'}, with $g : \mathbb{R}_+\times\mathbb{R}^3\rightarrow\mathbb{R}$ a function of space-time. Also assume that both $u$ and $g$ are smooth functions of $t \geqq 0$ and $x \in \mathbb{R}^{3}$ such that
\begin{equation}\label{2.3}
u, g \equiv 0 \quad \text { for } \quad|x| \geqq t+1
\end{equation}

In what follows we shall prove that, under suitable assumptions on $g$, we can recover some uniform decay rate for solutions to ~\eqref{eq:LKG}, ~\eqref{eq:IVP'}. Since our aim is to apply the results of this section to the nonlinear problem ~\eqref{eq:NKG} where the role of $g$ is played, typically, by a quadratic term, it is reasonable to assume that $g$, together with some of its derivatives, decays as a function of $t$ in the spatial $L^{2}$ norm. More precisely, using the notation introduced in section ~\eqref{section:introduction} we shall assume that $E_{k, N}(g)<\infty$ for some $0 \leqq k \leqq \frac{3}{2}$ and $N>0$, and prove that, with the same $k$, $|u|_{k, 0}<\infty$. For technical reasons we take $k$ to be $1+\varepsilon$, for some $0<\varepsilon \leqq \frac{1}{4}$, i.e.,
\begin{equation*}
    E_{1+\varepsilon, N}(g)<\infty   
\end{equation*}
for sufficiently large $N$. We prove:

\begin{proposition}\label{prop:1} There exists a constant $C>0$ such that the solution $u$ of ~\eqref{eq:LKG}, ~\eqref{eq:IVP'},  satisfies the inequality
\begin{equation}\label{eq:2.5}
|u(t, x)| \leqq C(1+t)^{-5 / 4}\left(\|u\|_{5}+E_{1+\varepsilon, 3}(g)\right)
\end{equation}

for any $t \geqq 0,|x| \leqq t+1$.
\end{proposition}
\begin{proof}
Performing the change of variable $t^{\prime}=t+2$ we may assume that the initial condition ~\eqref{eq:IVP'} is posed at $t'=2$ and that $u, g \equiv 0$ for $t' \geqq 2,|x| \geqq t'-1$. The crucial idea in our proof of ~\eqref{eq:2.5} is the use of pseudospherical coordinates in the interior of the light cone centered at $t'=0, x=0$. 

We introduce a change of coordinates from Cartesian coordinates \((t', x_1, x_2, x_3)\) to pseudo-spherical coordinates \((\rho, \theta, \phi, \psi)\). This choice is motivated by representing a point in \( \mathbb{R}^4 \) as \( \rho \omega \), where \( \rho \geq 0 \) is a radial component and \( w \) lies on the three-dimensional unit hyperboloid \( \mathbb{H} := \{(\omega_0, \omega_1, \omega_2, \omega_3) \in \mathbb{R}^4 \mid \omega_0^2 - \omega_1^2 - \omega_2^2 - \omega_3^2 = 1\} \). Using the identity $\cosh^2\theta - \sinh^2\theta = 1$, we define $\theta$ by setting \( \omega_0 = \cosh \theta \), \( \sqrt{\omega_1^2 + \omega_2^2 + \omega_3^2} = \sinh \theta \). Then we transform \( (\omega_1, \omega_2, \omega_3) \) into spherical coordinates on the two-dimensional sphere of radius \( \sinh \theta \) centered at the origin, with polar angle $\phi$ and azimuthal angle $\psi$. Thus we obtain the following change of coordinates:
\begin{align*}
    t' &= \rho \cosh \theta, \\
    x_1 &= \rho \sinh \theta \sin \phi \cos \psi, \\
    x_2 &= \rho \sinh \theta \sin \phi \sin \psi, \\
    x_3 &= \rho \sinh \theta \cos \phi.
\end{align*}
From now on, we use $t$ instead of $t'$ for simplicity. Throughout this proof, we will make use of multiple lemmas that will be proven at the end. 

\begin{lemma}\label{helper-lemma-1}
In these coordinates, the equation ~\eqref{eq:LKG} takes the form:
\[ \partial_\rho^2 u + \frac{3}{\rho} \partial_\rho u + u = \frac{1}{\rho^2} \Delta_H u + g, \]
with \(\Delta_H\) the Laplace-Beltrami operator on \(\mathbb{H}\) for the metric induced by the space-time metric \(ds^2 = dt^2 - dx_1^2 - dx_2^2 - dx_3^2\):
\[ \Delta_H = \sum_{i=1}^3 L_i^2 - \Delta_S, \]
where \(\Delta_S = \Omega_{12}^2 + \Omega_{23}^2 + \Omega_{31}^2\) is the Laplace-Beltrami operator of the 2-dimensional Euclidean unit sphere, $\mathbb{S}^2 = \{x \in \mathbb{R}^3 \mid \|x\| = 1\}$. 
\end{lemma}

\begin{lemma}\label{helper-lemma-2}
Introducing $v=\rho^{3 / 2} u$ we also have
\begin{equation}\label{eq:2.1''}
\partial_{\rho}^2v+v=h=h_{1}+h_{2}
\end{equation}
with
$$
h_{1}=\rho^{-1 / 2}\left(\Delta_{H} u+\frac{3}{4} u\right), \quad h_{2}=\rho^{3 / 2} g
$$
\end{lemma}
We define the cone $\Sigma = \{(t,x)\in\mathbb{R}_+\times\mathbb{R}^3 \mid \|x\| \leqq \frac{1}{2} t\}$, and the 3-dimensional hypersurface $\sigma = \Sigma \cap \mathbb{H}$. We can picture a section of these sets in the figure below, where the spacial domain is simplified to a 1-dimensional $x$-axis. The cone $\Sigma$ is represented in grey and the region where $u,g \neq 0$ is highlighted in yellow. 

\begin{figure}[H]
\centering
\begin{tikzpicture}
\begin{axis}[
    axis lines=middle,
    xlabel={$x$},
    ylabel={$t$},
    xmin=-3, xmax=3,
    ymin=0, ymax=5,
    xtick=\empty,
    ytick={1,2},
    xlabel style={at={(ticklabel* cs:1)},anchor=north west},
    ylabel style={at={(axis description cs:0.5,1)},anchor=south},
    clip=true,
    yticklabel style={font=\small, fill=white, inner sep=1pt},
    axis on top=true,
]

\addplot [domain=-4:4, samples=100, thick, dashed, name path=hyperboloid] {sqrt(x^2 + 1)};
\node at (axis cs:2.6,2.5) {$\mathbb{H}$}; 

\addplot [name path=SigmaTop, domain=-3:3, samples=100, draw=none] {5};
\addplot [name path=SigmaBottom, domain=-3:3, samples=100, draw=none] {2*abs(x)};
\addplot [gray!30, opacity=0.5] fill between[of=SigmaBottom and SigmaTop];
\node at (axis cs:1.8,4.6) {$\Sigma$}; 

\addplot [domain=-0.6:0.6, samples=100, thick, red] {sqrt(x^2 + 1)} node[pos=0.5, below right] {$\sigma$};

\addplot [name path=A, domain=-1:1, samples=100, thick] {2}; 
\addplot [name path=B, domain=-3:-1, samples=100, thick] {abs(x)+1}; 
\addplot [name path=C, domain=1:3, samples=100, thick] {abs(x)+1}; 

\addplot [name path=Atop, draw=none] coordinates {(-1,6) (1,6)}; 
\addplot [yellow, opacity=0.5] fill between[of=A and Atop]; 
\addplot [name path=Btop, draw=none] coordinates {(-3,6) (-1,6)};
\addplot [yellow, opacity=0.5] fill between[of=B and Btop];
\addplot [name path=TopRight, draw=none] coordinates {(1,6) (3,6)}; 
\addplot [yellow, opacity=0.5] fill between[of=C and TopRight];
\node at (axis cs:-2.3,4.6) {$u,g\not\equiv0$};

\end{axis}
\end{tikzpicture}
\end{figure}

\begin{lemma}\label{helper-lemma-3}
For each $\omega \in \sigma$, the solution of the inhomogeneous linear second order differential equation ~\eqref{eq:2.1''}, ~\eqref{eq:IVP'},  is given by:
\begin{align*}
v(\rho \omega) = v_0(\rho \omega) + \int_{\rho_0}^\rho \sin(\rho-\lambda) h(\lambda \omega) d\lambda
\end{align*}
where
\begin{align*}
    v_0(\rho \omega) = \cos(\rho-\rho_0)v(\rho_0 \omega) + \sin(\rho-\rho_0) \partial_\rho v(\rho_0 \omega)
\end{align*}
and $\rho_{0}$ is the $\rho$-coordinate of the corresponding point on $t=2$. 
\end{lemma}

Clearly, we have,
\begin{equation}\label{eq:2.8}
\left|v_{0}(\rho \omega)\right|^2 \leqq C \left(\left\|u_{0}\right\|_{L^2}^2 + \left\|u_{1}\right\|_{L^2}^2 + \sum_{i=1}^3 \left\|\partial_i u_{0}\right\|_{L^2}^2 \right) := CD^2
\end{equation}
for every $\rho \geqq \rho_{0}, \omega \in \sigma$.
\begin{lemma}\label{helper-lemma-4}
On the other hand, since $u, g \equiv 0$ for $|x| \geqq t-1, t \geqq 2$, and for $x \in \mathbb{R}^{3}, t>2$, we have, for any $\omega \in \mathbb{H} \backslash \sigma$,
\begin{align*}
v(\rho \omega)=\int_{1}^{\rho} \sin (\lambda-\rho) h(\lambda \omega) d \lambda
\end{align*}
\end{lemma}
Now, using the Cauchy-Schwartz inequality, we have, for every $\omega\in\sigma$, 
\begin{align*}
|v(\rho \omega)|^{2} &\leqq D^{2}+\int_{1}^{\rho} \lambda^{-1-\varepsilon}\lambda^{1+\varepsilon}|h(\lambda \omega)|^{2} d \lambda \\
&\leqq C\left( D^{2}+\int_{1}^{\rho} \lambda^{1+\varepsilon}|h(\lambda \omega)|^{2} d \lambda \right)
\end{align*}
with $C:=\int_1^\infty\lambda^{-1-\varepsilon}d\lambda<\infty$, and $\varepsilon>0$ a fixed, small, positive constant. Integrating over $\sigma$, and changing $C$, we obtain
\begin{align*}
\int_{\sigma}|v(\rho \omega)|^{2} d S_{\omega} \leqq C\left(D^{2}+\int_{\sigma}\int_{1}^{\rho}   \lambda^{1+\varepsilon}|h(\lambda \omega)|^{2} d \lambda d S_{\omega}\right)
\end{align*}
We recall now the definitions of $h$ in ~\eqref{eq:2.1''} and the form of the volume element in $\mathbb{R}^{4}$ expressed in pseudospherical coordinates, $dtd x=\lambda^{3} d \lambda d S_{\omega}$, with $d S_{\omega}$ the area element on $\mathbb{H}$. Thus, since $u, g$ vanish for $|x| \leqq t-1, t \geqq 2$, and for $t<2$, we derive
\begin{align*}
    \int_{\sigma}&|v(\rho \omega)|^{2} d S_{\omega} 
    \leqq C\left(D^{2}+\int_{\sigma}\int_{1}^{\rho} \lambda^{-2+\varepsilon}|h(\lambda \omega)|^{2} \lambda^3 d \lambda d S_{\omega}\right) \\
    &\leqq C\left(D^{2}+\int_{s=1}^{\frac{1+\rho^2}{2}}\int_{|y|\leqq s-1} (s^2-|y|^2)^\frac{-2+\varepsilon}{2}\left|(s^2-|y|^2)^{-\frac{1}{2}}\left(\Delta_H u(s,y) + u(s,y)\right)\right|^{2}dyds \right. \\
    & \quad \quad + \left. \int_{s=1}^{\frac{1+\rho^2}{2}}\int_{|y|\leqq s-1} (s^2-|y|^2)^\frac{-2+\varepsilon}{2}\left|(s^2-|y|^2)^\frac{3}{2}g(s,y)\right|^{2}dyds\right) \\
    &\leqq C\left(D^{2}+\int_{s=1}^{\frac{1+\rho^2}{2}}\int_{|y|\leqq s-1} (s^2-|y|^2)^\frac{-3+\varepsilon}{2}\left(\left|\Delta_H u(s,y)\right|^2 + \left|u(s,y)\right|^2 \right) dyds \right. \\
    & \quad \quad + \left. \int_{s=1}^{\frac{1+\rho^2}{2}}\int_{|y|\leqq s-1} (s^2-|y|^2)^\frac{1+\varepsilon}{2}\left|g(s,y)\right|^{2}dyds\right) 
\end{align*}
Now, recalling the definition of $\|u\|_2$ and absorbing $D^2$ into it,
\begin{align*}
    D^2 & + \int_{s=1}^{\frac{1+\rho^2}{2}} \int_{|y|\leqq s-1} (s^2-|y|^2)^\frac{-3+\varepsilon}{2}\left(\left|\Delta_H u(s,y)\right|^2 + \left|u(s,y)\right|^2 \right) dyds \\
    &\leqq D^2 + \int_{s=1}^{\frac{1+\rho^2}{2}}\int_{|y|\leqq s-1} (s-|y|)^\frac{-3+\varepsilon}{2} (s+|y|)^\frac{-3+\varepsilon}{2} \left(\left|\Delta_H u(s,y)\right|^2 + \left|u(s,y)\right|^2 \right) dyds \\
    &\leqq D^2 + \int_{s=1}^{\frac{1+\rho^2}{2}} s^\frac{-3+\varepsilon}{2} \int_{|y|\leqq s-1} \left(\left|\Delta_H u(s,y)\right|^2 + \left|u(s,y)\right|^2 \right) dyds \quad \leqq \quad C \|u\|_2^2
\end{align*}
With $C:=\int_1^\infty s^\frac{-3+\varepsilon}{2} ds<\infty$. 
And similarly, recalling the definition of $E_{1+\varepsilon, 0}(g)$,
\begin{align*}
    \int_{s=1}^{\frac{1+\rho^2}{2}}\int_{|y|\leqq s-1} & (s^2-|y|^2)^\frac{1+\varepsilon}{2}\left|g(s,y)\right|^{2}dyds \\
    &\leqq \int_{s=1}^{\frac{1+\rho^2}{2}} (1+s)^\frac{-3-3\varepsilon}{2} \left((1+s)^{1+\varepsilon}\right)^2 \int_{|y|\leqq s-1}  \left|g(s,y)\right|^{2}dyds \\
    &\leqq CE_{1+\varepsilon, 0}^2(g) 
\end{align*}
with $C:=\int_1^\infty (1+s)^\frac{-3-3\varepsilon}{2} ds <\infty$.
Thus we obtain
\begin{align}\label{eq:13}
    \int_{\sigma}|v(\rho \omega)|^{2} d S_{\omega} \leqq C\left(\|u\|_2^2 + E_{1+\varepsilon, 0}^2(g)\right)
\end{align}
Now apply the classical Sobolev inequality ~\eqref{lemma:injection} on the unit hyperboloid $\mathbb{H}$, i.e., given $v(\omega)$ defined, in a neighborhood $\sigma^{\prime}$ of $\sigma$ in $\mathbb{H}$,
\begin{equation*}
\sup _{\omega \in \sigma}|v(\omega)| \leqq C\|v\|_{W^{s, 2}\left(\sigma^{\prime}\right)}, \quad s > \frac{3}{2},
\end{equation*}
with $W^{s, 2}\left(\sigma^{\prime}\right)$ the usual Sobolev spaces on $\sigma^{\prime} \subset \mathbb{H}$. Since the angular derivatives on $\mathbb{H}$ can be expressed with the help of the operators $\{\Omega_{a b}\mid0 \leq a<b \leq 3\}$, the norms $\|\cdot\|_{W^{s, 2}\left(\sigma^{\prime}\right)}$ are equivalent to 
\begin{align*}
    \left(\sum_{|\beta| \leqq s}\left\|\Omega^{\beta} v\right\|_{L^{2}\left(\sigma^{\prime}\right)}^{2}\right)^{1/2},
\end{align*}
with  $\beta=(\beta_{ab})_{0 \leq a<b \leq 3}$ being 6-indices of size $|\beta|=\sum \beta_{ab}$, and 
\begin{align*}
    \Omega^{\beta}=\prod \Omega_{a b}^{\beta_{ab}}.
\end{align*}
Thus we have, for every $\rho \geqq 1$,
\begin{align*}
\sup _{\omega \in \sigma}|v(\rho \omega)| \leqq C\left(\sum_{|\beta| \leqq 2} \int_{\sigma^{\prime}}\left|\Omega^{\beta} v(\rho \omega)^{2}\right| d S_{\omega}\right)^{1 / 2}
\end{align*}
Therefore, making use of the commutation properties of the $\Omega_{a b}$ 's with $\rho$ and $\square_{1}$, we derive, from the inequality ~\eqref{eq:13} applied on a slightly larger domain $\sigma^{\prime}>\sigma$ in $\mathbb{H}$,
\begin{equation}\label{eq:sup}
\sup _{\omega \in \sigma}|v(\rho \omega)| \leqq C\left(\|u\|_{4}+E_{1+\varepsilon, 2}(g)\right)
\end{equation}
for any $\rho \geqq 1$, with $C$ independent of $\rho$. i.e.,
\begin{equation*}
|v(\rho \omega)| \leqq C\left(\|u\|_{4}+E_{1+\varepsilon, 2}(g)\right)
\end{equation*}
for any $\rho\geqq1$, $\omega \in \sigma$. Finally, since $v=\rho^{3 / 2} u$, and recalling that $t\geq2|x|$, we obtain
\begin{align*}
|u(t, x)| &\leqq C\left(t^{2}-|x|^{2}\right)^{-3 / 4}\left(\|u\|_{4}+E_{1+\varepsilon, 2}(g)\right) \\
&\leqq C\left(t-|x|\right)^{-3 / 4}\left(t+|x|\right)^{-3 / 4}\left(\|u\|_{4}+E_{1+\varepsilon, 2}(g)\right) \\
&\leqq C\left(1+t\right)^{-6 / 4}\left(\|u\|_{4}+E_{1+\varepsilon, 2}(g)\right) \\
\end{align*}
for all $(t, x) \in \Sigma, t \geqq 2$, which proves ~\eqref{eq:2.5} in that region.

To prove ~\eqref{eq:2.5} outside of $\Sigma$ we have to rely on a different form of the Sobolev inequality than that used above. With this in mind, write, for any $\omega\in\mathbb{H}$, $\omega=\left(\omega_{0}, \omega^{\prime}\right)$ with $\omega_{0}=\cosh \theta, \omega^{\prime}=\left(\omega_{1}, \omega_{2}, \omega_{3}\right)$ $=(\sinh \theta) \xi$ with $\xi\in\mathbb{S}^2$ and $\theta\in\mathbb{R}_+$. Given $\omega \in \mathbb{H}$, we define $H_{\omega}=\{\tilde{\omega} \in \mathbb{H} \mid \tilde{\theta} \geqq \theta\}$. Proceeding precisely as in the derivation of ~\eqref{eq:13} we infer, using lemma ~\eqref{helper-lemma-4}
\begin{equation}\label{eq:2.10'}
\int_{H_{\omega}}|v(\rho \tilde{\omega})|^{2} d S_{\tilde{\omega}} \leqq C\left(\|u\|_{2}^{2}+E_{1+\varepsilon, 0}^2(g)\right)
\end{equation}
for any $\omega \in H \backslash \sigma$. Now we make use of the following:
\begin{lemma}\label{helper-lemma-5}
Let $v=v(\omega)$ be a smooth, compactly supported function in $\mathbb{H}$. There exists a constant $C>0$ such that, for any $\omega \in \mathbb{H}, \omega=(\cosh \theta,(\sinh \theta) \xi)$, with $|\xi|=1$ and $\theta>0$,
\begin{align*}  
    |v(\omega)| \leqq C \frac{1}{\sinh \theta}\|v\|_{W^{3,2}\left(H_{\omega}\right)}
\end{align*}
\end{lemma}
Applying this lemma as in the proof of ~\eqref{eq:sup} to $v(\rho \omega)$, which is compactly supported in $\mathbb{H}$ for any fixed $\rho$, we infer from ~\eqref{eq:2.10'}
\begin{equation*}
|v(\rho \omega)| \leqq C \frac{1}{\sinh \theta}\left(\|u\|_{5}+E_{1+e, 3}(g)\right)
\end{equation*}
for any $\omega \in \mathbb{H} \backslash \sigma$, with $\omega=(\cosh \theta,(\sinh \theta) \xi)$, and any $\rho \geqq 1$.
Finally, since $v=\rho^{3 / 2} u$, and $|x|=\rho \sinh \theta$, we obtain
\begin{align*}
    |u(t, x)| &\leqq C \frac{1}{\rho^{3 / 2} \sinh \theta}\left(\|u\|_{5}+E_{1+\varepsilon, 3}(g)\right) \\
    &\leqq C (t^2-|x|^2)^{-\frac{1}{4}}|x|^{-1}\left(\|u\|_{5}+E_{1+\varepsilon, 3}(g)\right) \\
    &\leqq C (t-|x|)^{-\frac{1}{4}}(t+|x|)^{-\frac{1}{4}}|x|^{-1}\left(\|u\|_{5}+E_{1+\varepsilon, 3}(g)\right) \\
    &\leqq C (t+|x|)^{-\frac{1}{4}}(t-1)^{-1}\left(\|u\|_{5}+E_{1+\varepsilon, 3}(g)\right) \\
    &\leqq C (1+t)^{-\frac{5}{4}}\left(\|u\|_{5}+E_{1+\varepsilon, 3}(g)\right) \\
\end{align*}
for any $(t, x)$ such that $t-1 \geqq|x| \geqq \frac{1}{2} t, t \geqq 2$ which implies, in particular, ~\eqref{eq:2.5} and thus ends the proof of Proposition ~\ref{prop:1}.
\end{proof}

\begin{proof}[Proof of Lemma~\eqref{helper-lemma-1}]
We compute:
\begin{align*}
    \partial_\rho = \partial_\rho t \cdot \partial_t + \sum_{i=1}^3 \partial_\rho x_i \cdot \partial_{i} = \frac{t}{\rho} \partial_t + \sum_{i=1}^3 \frac{x_i}{\rho} \partial_{i} = \frac{1}{\rho} \left(t \partial_t + \sum_{i=1}^3 x_i \partial_i \right), 
\end{align*}
\begin{align*}
& \partial_t\left[\frac{t}{\rho}\right] = \partial_t\left[\frac{t}{\sqrt{t^2-|x|^2}}\right] = -\frac{|x|^2}{(t^2 - |x|^2)^{\frac{3}{2}}} =  -\frac{|x|^2}{\rho^3},\\
& \partial_t\left[\frac{x_i}{\rho}\right] = \partial_t\left[\frac{x_i}{\sqrt{t^2 - |x|^2}}\right] = -\frac{x_i t}{(t^2 - |x|^2)^{\frac{3}{2}}} = -\frac{x_i t}{\rho^3},\\
& \partial_i\left[\frac{t}{\rho}\right] = \partial_i\left[\frac{t}{\sqrt{t^2 - |x|^2}}\right] = \partial_i\left[\frac{t}{\sqrt{t^2 - \sum x_j^2}}\right] = \frac{t x_i}{\rho^3}, \\
& \partial_i\left[\frac{x_j}{\rho}\right] = \frac{\delta_{ij}\rho -x_j\partial_i\rho}{\rho^2} = \frac{\delta_{ij}}{\rho} + \frac{x_i x_j}{\rho^3}.
\end{align*}
\begin{align*}
\partial_\rho^2 &= \frac{1}{\rho^2} \left(t \partial_t + \sum_{i=1}^3 x_i \partial_i \right)^2 \\
&= \frac{1}{\rho^2} \left( t^2 \partial_t^2 + 2t \sum_{i=1}^3 x_i \partial_i \partial_t + \left(\sum_{i=1}^3 x_i \partial_i\right) \left(\sum_{j=1}^3 x_j \partial_j\right) \right) \\
&= \frac{1}{\rho^2} \left( t^2 \partial_t^2 + 2t \sum_{i=1}^3 x_i \partial_i \partial_t + \sum_{i,j=1}^3 x_i x_j \partial_i \partial_j \right).
\end{align*}
\begin{align*}
    \Omega_{ij} &= x_i\partial_j - x_j\partial_i,\\
    \Omega_{ij}^2 &= (x_i\partial_j - x_j\partial_i)(x_i\partial_j - x_j\partial_i) = x_i^2\partial_j^2 - x_i\partial_i - x_ix_j\partial_i\partial_j - x_j\partial_j - x_jx_i\partial_j\partial_i + x_j^2\partial_i^2 \\
    &= x_i^2\partial_j^2 + x_j^2\partial_i^2 - 2x_ix_j\partial_i\partial_j - x_i\partial_i - x_j\partial_j \\
    \Omega_{0i} &= t\partial_i + x_i\partial_t, \\
    \Omega_{0i}^2 &= (t\partial_i + x_i\partial_t)(t\partial_i + x_i\partial_t) = t^2\partial_i^2 + t \partial_t + tx_i\partial_i\partial_t + x_i\partial_i + tx_i\partial_t\partial_i + x_i^2\partial_t^2 \\
    &= t^2\partial_i^2 + x_i^2\partial_t^2 + 2tx_i\partial_t\partial_i + t\partial_t + x_i\partial_i.
\end{align*}
\begin{align*}
\Delta_H &= \left(\sum_{i=1}^3 \Omega_{0i}^2\right) - \Delta_S = \left(\sum_{i=1}^3 \Omega_{0i}^2\right) - \left(\Omega_{12}^2 + \Omega_{23}^2 + \Omega_{31}^2\right) \\
&= t^2 \sum_{i=1}^3 \partial_i^2 + \sum_{i=1}^3 x_i^2 \partial_t^2 + 2t \sum_{i=1}^3 x_i \partial_t \partial_i + 3t \partial_t + \sum_{i=1}^3 x_i \partial_i - x_1^2 \partial_2^2 - x_2^2 \partial_1^2 + 2x_1 x_2 \partial_1 \partial_2 \\
& \quad + x_1 \partial_1 + x_2 \partial_2 - x_2^2 \partial_3^2 - x_3^2 \partial_2^2 + 2x_2 x_3 \partial_2 \partial_3 + x_2 \partial_2 + x_3 \partial_3 - x_3^2 \partial_1^2 - x_1^2 \partial_3^2 \\
& \quad \quad + 2x_3 x_1 \partial_3 \partial_1 + x_3 \partial_3 + x_1 \partial_1 \\
&= t^2 \sum_{i=1}^3 \partial_i^2 + |x|^2 \partial_t^2 + 2t \sum_{i=1}^3 x_i \partial_t \partial_i + 3t \partial_t + \sum_{i=1}^3 x_i \partial_i - \sum_{\substack{i,j=1 \\ i \neq j}}^3 x_i^2 \partial_j^2 \\
& \quad \quad + \sum_{\substack{i,j=1 \\ i \neq j}}^3 x_i x_j \partial_i \partial_j + 2 \sum_{i=1}^3 x_i \partial_i \\
&= t^2 \sum_{i=1}^3 \partial_i^2 + |x|^2 \partial_t^2 + 2t \sum_{i=1}^3 x_i \partial_t \partial_i + 3t \partial_t + 3 \sum_{i=1}^3 x_i \partial_i - \sum_{i,j=1}^3 x_i^2 \partial_j^2 \\
& \quad \quad + \sum_{i=1}^3 x_i^2 \partial_i^2 + \sum_{i,j=1}^3 x_i x_j \partial_i \partial_j - \sum_{i=1}^3 x_i^2 \partial_i^2 \\
&= t^2 \Delta + |x|^2 \partial_t^2 + 2t \sum_{i=1}^3 x_i \partial_t \partial_i + 3t \partial_t + 3 \sum_{i=1}^3 x_i \partial_i - \sum_{i,j=1}^3 x_i^2 \partial_j^2 + \sum_{i,j=1}^3 x_i x_j \partial_i \partial_j.
\end{align*}
Putting all together, and using ~\eqref{eq:LKG}, we get:
\begin{align*}
\frac{1}{\rho^2} \Delta_H u + & g - u = \frac{1}{\rho^2} \Delta_H u + \partial_t^2 u - \Delta u \\
&= \frac{1}{\rho^2} \left( (t^2 - \rho^2) \Delta u + (|x|^2 + \rho^2) \partial_t^2 u + 2t \sum_{i=1}^3 x_i \partial_t \partial_i u + 3t \partial_t u + 3 \sum_{i=1}^3 x_i \partial_i u \right. \\
&\quad \left. - \sum_{i,j=1}^3 x_i^2 \partial_j^2 u + \sum_{i,j=1}^3 x_i x_j \partial_i \partial_j u \right) \\
&= \frac{1}{\rho^2} \left( |x|^2 \Delta u + t^2 \partial_t^2 u + 2t \sum_{i=1}^3 x_i \partial_t \partial_i u + 3t \partial_t u + 3 \sum_{i=1}^3 x_i \partial_i u \right. \\
&\quad \left. - \sum_{i,j=1}^3 x_i^2 \partial_j^2 u + \sum_{i,j=1}^3 x_i x_j \partial_i \partial_j u \right) \\
&= \frac{1}{\rho^2} \left( t^2 \partial_t^2 u + 2t \sum_{i=1}^3 x_i \partial_t \partial_i u + \sum_{i,j=1}^3 x_i x_j \partial_i \partial_j u + 3t \partial_t u + 3 \sum_{i=1}^3 x_i \partial_i u \right) \\
&= \partial_{\rho}^2 u + \frac{3}{\rho} \partial_{\rho} u 
\end{align*}
\end{proof}

\begin{proof}[Proof of Lemma ~\eqref{helper-lemma-2}]
\begin{align*}
\partial_\rho v &= \partial_\rho \left[\rho^{3/2} u \right] = \frac{3}{2} \rho^{1/2} u + \rho^{3/2} \partial_\rho u, \\
\partial_\rho^2 v &= \frac{3}{2} \frac{1}{2} \rho^{-1/2} u + \frac{3}{2} \rho^{1/2} \partial_\rho u + \frac{3}{2} \rho^{1/2} \partial_\rho u + \rho^{3/2} \partial_\rho^2 u \\
&= \frac{3}{4} \rho^{-1/2} u + 3 \rho^{1/2} \partial_\rho u + \rho^{3/2} \partial_\rho^2 u \\
&= \frac{3}{4} \rho^{-1/2} u + 3 \rho^{1/2} \partial_\rho u + \rho^{3/2} \left(\frac{1}{\rho^2} \Delta_H u + g - \frac{3}{\rho} \partial_\rho u - u\right) \\
&= \frac{3}{4} \rho^{-1/2} u + 3 \rho^{1/2} \partial_\rho u + \rho^{-1/2} \Delta_H u + \rho^{3/2} g - 3 \rho^{1/2} \partial_\rho u - \rho^{3/2} u \\
&= \frac{3}{4} \rho^{-1/2} u + \rho^{-1/2} \Delta_H u + \rho^{3/2} g - \rho^{3/2} u \\
&=  \rho^{-1/2}\left[ \frac{3}{4} u + \Delta_H u \right] + \rho^{3/2} g - v. 
\end{align*}

\end{proof}

\begin{proof}[Proof of Lemma~\ref{helper-lemma-3}]
For a fixed $w \in \sigma$, set $V(\rho) := v(\rho w)$ and $H(\rho) := h(\rho w)$ for all $\rho\geqq0$. We have:
\begin{align*}
    H = \left(\partial_\rho^2 + 1\right) V = -\left(\left(\frac{1}{i}\partial_\rho\right)^2-1\right) V = -\left(\frac{1}{i}\partial_\rho-1\right) \left(\frac{1}{i}\partial_\rho+1\right) V
\end{align*}
Then set $W:=\left(\frac{1}{i}\partial_\rho+1\right) V$ and let us solve the following first oder differential equation for $W$:
\begin{align*}
    H = -\left(\frac{1}{i}\partial_\rho-1\right) W \Leftrightarrow -iH = \partial_\rho W - iW
\end{align*}
Note that
\begin{align*}
    \partial_\rho \left[ e^{-i\rho} W(\rho) \right] = -i e^{-i\rho} W(\rho) + e^{-i\rho} \partial_\rho W(\rho) = e^{-i\rho}(-i W(\rho) + \partial_\rho W(\rho)) = -ie^{-i\rho} H(\rho)
\end{align*}
Let $\rho_0$ be the $\rho$-coordinate of the point corresponding to $w$, at $t=2$. Integrating from $\rho_0$ to $\rho\geqq\rho_0$ yields:
\begin{align*}
    & e^{-i\rho} W(\rho) = e^{-i\rho_0} W(\rho_0) - i \int_{\rho_0}^\rho e^{-is} H(s) ds \\
    & W(\rho) = e^{i(\rho-\rho_0)} W(\rho_0) - i \int_{\rho_0}^\rho e^{i(\rho-s)} H(s) ds := K(\rho) 
\end{align*}
So we have
\begin{align*}
    \left(\frac{1}{i}\partial_\rho+1\right) V = K \Leftrightarrow \partial_\rho V + iV = iK
\end{align*}
And as earlier, note that
\begin{align*}
    \partial_\rho\left[ e^{i\rho} V(\rho)\right] = ie^{i\rho}V(\rho) + e^{i\rho}\partial_\rho V(\rho) = e^{i\rho} (iV(\rho)+\partial_\rho V(\rho)) = ie^{i\rho}K(\rho)
\end{align*}
Integrate from $\rho_0$ to $\rho\geqq\rho_0$:
\begin{align*}
    & e^{i\rho} V(\rho) = e^{i\rho_0} V(\rho_0) + i \int_{\rho_0}^\rho e^{i\lambda} K(\lambda) d\lambda \\
    & V(\rho) = e^{-i(\rho-\rho_0)} V(\rho_0) + i \int_{\rho_0}^\rho e^{-i(\rho-\lambda)} K(\lambda) d\lambda 
\end{align*}
We compute:
\begin{align*}
    \int_{\rho_0}^\rho e^{-i(\rho-\lambda)} & K(\lambda) d\lambda = \int_{\rho_0}^\rho e^{-i(\rho-\lambda)} \left( e^{i(\lambda-\rho_0)} W(\rho_0) - i \int_{\rho_0}^\lambda e^{i(\lambda-s)} H(s) ds \right) d\lambda \\
    &= e^{-i(\rho+\rho_0)}W(\rho_0)\int_{\rho_0}^\rho e^{2i\lambda}d\lambda - i \int_{\rho_0}^\rho e^{2i\lambda} \left( \int_{\rho_0}^\lambda e^{-i(\rho+s)} H(s) ds \right) d\lambda \\
    &= e^{-i(\rho+\rho_0)}W(\rho_0)\left[ \frac{e^{2i\lambda}}{2i}\right]_{\lambda=\rho_0}^\rho - i \int_{\rho_0}^\rho \partial_\lambda \left[ \frac{e^{2i\lambda}}{2i}\right] \left( \int_{\rho_0}^\lambda e^{-i(\rho+s)} H(s) ds \right) d\lambda \\
    &= \frac{W(\rho_0)}{2i} \left[ e^{i(\rho-\rho_0)} - e^{-i(\rho-\rho_0)}\right] - i \left[ \frac{e^{2i\lambda}}{2i} \int_{\rho_0}^\lambda e^{-i(\rho+s)} H(s) ds\right]_{\lambda=\rho_0}^\rho \\
    & \quad + i \int_{\rho_0}^\rho \frac{e^{2i\lambda}}{2i} e^{-i(\rho+\lambda)} H(\lambda) d\lambda \\
    &= W(\rho_0)\sin(\rho-\rho_0) -  i \int_{\rho_0}^\rho \frac{e^{-i(\rho-s)}}{2i} H(s) ds + i \int_{\rho_0}^\rho \frac{e^{-i(\rho-\lambda)}}{2i} H(\lambda) d\lambda \\
    &= V(\rho_0)\sin(\rho-\rho_0) + \frac{1}{i}\partial_\rho V(\rho_0)\sin(\rho-\rho_0) - i \int_{\rho_0}^\rho \sin(\rho-\lambda) H(\lambda) d\lambda
\end{align*}
Thus, we get:
\begin{align*}
    V(\rho) &= e^{-i(\rho-\rho_0)} V(\rho_0) + iV(\rho_0)\sin(\rho-\rho_0) + \partial_\rho V(\rho_0)\sin(\rho-\rho_0) + \int_{\rho_0}^\rho \sin(\rho-\lambda) H(\lambda) d\lambda \\
    &= \cos(\rho-\rho_0)V(\rho_0) + \sin(\rho-\rho_0) \partial_\rho V(\rho_0) + \int_{\rho_0}^\rho \sin(\rho-\lambda) H(\lambda) d\lambda
\end{align*}
\end{proof}

\begin{proof}[Proof of Lemma ~\eqref{helper-lemma-4}]
    The proof follows the exact same arguments as the previous one, but integrating from $1$ to $\rho\geqq1$ instead of integrating from $\rho_0$ to $\rho\geqq\rho_0$, and recalling that $V(1) = \partial_\rho V(1) = 0$. 
\end{proof}

\begin{proof}[Proof of Lemma ~\eqref{helper-lemma-5}]
For any $\omega=(\cosh \theta,(\sinh \theta) \xi)$ and $v = v(\omega)$, we set $f(\theta, \xi) = v(\omega) = v(\cosh \theta,(\sinh \theta) \xi)$. Then we write (using the fundamental theorem of calculus and the fact that $f$ has compact support):

\begin{align*}
    f^2(\theta,\xi) = 2 \int_{\theta}^{\infty} f(\tilde{\theta}, \xi) \cdot \partial_{\theta}f(\tilde{\theta}, \xi) d \tilde{\theta}
\end{align*}
Hence, for $\theta > 0$, and since $\sinh$ is non-decreasing on $\mathbb{R}_+\backslash\{0\}$,
\begin{align*}
    |f(\theta,\xi)|^{2} \leqq \frac{2}{\sinh ^{2} \theta} \int_{\theta}^{\infty}\left|f(\tilde{\theta}, \xi) \cdot \partial_{\theta}f(\tilde{\theta}, \xi) \right| \sinh ^{2} \tilde{\theta} d \tilde{\theta}
\end{align*}
Integrating on the unit sphere $\mathbb{S}^2 = \{\xi \in \mathbb{R}^3 \mid \|\xi\|=1\}$, we get:
\begin{align*}
\|f(\theta,\;\cdot\;)\|^2_{L^{2}\left(\mathbb{S}^2\right)} = \int_{\|\xi\|=1}|f(\theta,\xi)|^{2} d S_{\xi} \leqq \frac{2}{\sinh ^{2} \theta} \int_{\|\xi\|=1} \int_{\theta}^{\infty}\left|f(\tilde{\theta}, \xi) \cdot \partial_{\theta}f(\tilde{\theta}, \xi) \right| \sinh ^{2} \tilde{\theta} d \tilde{\theta} d S_{\xi}
\end{align*}
Then, applying the Cauchy-Schwartz inequality yields
\begin{align}
\|f(\theta, \;\cdot\;)\|^2_{L^{2}(\mathbb{S}^2)} &\leqq \frac{2}{\sinh^{2} \theta}\|f\|_{L^{2}(H_{\omega})} \cdot \|\partial_{\theta}f\|_{L^2(H_{\omega})} \label{eq:inequality1} \\
&\leqq \frac{2}{\sinh^{2} \theta}\|f\|_{L^{2}(H_{\omega})} \cdot \|f\|_{W^{1,2}(H_{\omega})} \label{eq:inequality2}
\end{align}
Now, notice that
\begin{align*}
|f(\theta, \xi)|^2 \leqq \sup_{|\xi|=1} |f(\theta, \xi)|^2 = \|f(\theta, \;\cdot\;)\|^2_{L^\infty(\mathbb{S}^2)}
\end{align*}
and using lemma ~\ref{lemma:injection}, we have
\begin{align*}
    \|f(\theta, \;\cdot\;)\|^2_{L^\infty(S^2)} \leqq \|f(\theta, \;\cdot\;)\|^2_{W^{2,2}(\mathbb{S}^2)} = \sum_{|\alpha| \leq 2} \|D^{\alpha}f(\theta, \;\cdot\;)\|^2_{L^2(\mathbb{S}^2)}
\end{align*}
Now, using \eqref{eq:inequality2} with \(D^{\alpha}f(\theta, \;\cdot\;)\), we have
\begin{align*}
    \|D^{\alpha}f(\theta, \;\cdot\;)\|^2_{L^{2}(\mathbb{S}^2)} \leqq \frac{2}{\sinh^{2} \theta}\|D^{\alpha} f\|_{L^{2}(H_{\omega})} \cdot \|D^{\alpha} f\|_{W^{1,2}(H_{\omega})}
\end{align*}
thus, we get
\begin{align*}
|f(\theta, \xi)|^2 \leqq \frac{2}{\sinh^{2} \theta} \sum_{|\alpha| \leq 2} \|D^{\alpha} f\|_{L^{2}(H_{\omega})} \cdot \|D^{\alpha} f\|_{W^{1,2}(H_{\omega})} \leqq \frac{2}{\sinh^{2} \theta} \|f\|^2_{W^{3,2}(H_{\omega})}
\end{align*}
\end{proof}

We can get rid of the term $\|u\|_{5}$ on the right-hand side of ~\eqref{eq:2.5} by using the following straightforward consequence of the energy identity for ~\eqref{eq:LKG}:
\begin{proposition}\label{prop:3-1}
Let $u$ be a solution of ~\eqref{eq:LKG} satisfying ~\eqref{eq:IVP'}. Then for any $N \in \mathbb{N}$, and for any sufficiently small $\varepsilon>0$,
\begin{equation}
    \|u\|_{N} \leqq C(\Delta_N + E_{1+\varepsilon, N}(g)),
\end{equation}
\end{proposition}
\begin{proof}
Let $\alpha$ be a 10-index of size $|\alpha| \leqq N$. Applying $\Gamma^{\alpha}$ to both sides of ~\eqref{eq:LKG} and since $[\Gamma^{\alpha}, \Box_1] = 0$ (see lemma ~\eqref{lemma:GammaAlphaBox1}), we get 
$$\square_1 \Gamma^{\alpha} u = \Gamma^{\alpha} g.$$
We then multiply both sides of the equation by $\partial_t \Gamma^{\alpha} u$ and integrate over the spatial domain $\mathbb{R}^3$. By lemma ~\eqref{lemma:energy-helper}, the LHS is equal to $\frac{1}{2}\partial_t \mathcal{E}^2(\Gamma^{\alpha} u (t, \;\cdot\;))$. Using the Cauchy-Schwartz inequality on the RHS, we get:
\begin{align*}
\frac{1}{2}\partial_t \mathcal{E}^2(\Gamma^{\alpha} u (t, \;\cdot\;)) = \int_{\mathbb{R}^3} \partial_t \Gamma^{\alpha} u(t,x) \cdot \Gamma^{\alpha} g(t,x) \, dx &\leqq \|\partial_t \Gamma^{\alpha} u (t, \;\cdot\;)\|_{L^2} \cdot \|\Gamma^{\alpha} g(t, \;\cdot\;)\|_{L^2} \\
&\leqq \mathcal{E}(\Gamma^{\alpha} u(t, \;\cdot\;)) \cdot \|\Gamma^{\alpha} g(t, \;\cdot\;)\|_{L^2}
\end{align*}
Thus, we have 
\begin{align*}
    \frac{1}{2}\partial_t \mathcal{E}^2(\Gamma^{\alpha} u (t, \;\cdot\;)) =  \mathcal{E}(\Gamma^{\alpha} u (t, \;\cdot\;)) \cdot \partial_t \mathcal{E}(\Gamma^{\alpha} u (t, \;\cdot\;))\leqq \mathcal{E}(\Gamma^{\alpha} u(t, \;\cdot\;)) \cdot \|\Gamma^{\alpha} g(t, \;\cdot\;)\|_{L^2}
\end{align*}
i.e.,
\begin{align*}
    \partial_t \mathcal{E}(\Gamma^{\alpha} u(t, \;\cdot\;)) \leqq \|\Gamma^{\alpha} g(t, \;\cdot\;)\|_{L^2}
\end{align*}
Integrating with respect to time $s$, from $s=0$ to $s=t$, yields:
\begin{align*}
    \mathcal{E}(\Gamma^{\alpha} u(t, \;\cdot\;)) \leqq \mathcal{E}(\Gamma^{\alpha} u(0, \;\cdot\;)) + \int_0^t \|\Gamma^{\alpha} g(s, \;\cdot\;)\|_{L^2}\, ds.
\end{align*}
We then take the sum over all 10-indices $\alpha$ of size $|\alpha| \leqq N$ to get:
\begin{align*}
    \|u(t, \;\cdot\;)\|_{\Gamma, N} &\leqq \sum_{|\alpha| \leqq N} \mathcal{E}(\Gamma^{\alpha} u (0, \;\cdot\;)) + \int_0^t \sum_{|\alpha| \leqq N} \|\Gamma^{\alpha} g(s, \;\cdot\;)\|_{L^2}\, ds \\
    &\leqq \Delta_N +  \int_0^t \|g(s, \;\cdot\;)\|_{\Gamma, N}\, ds
\end{align*}
Now, let $\varepsilon>0$ be a small parameter and we have
\begin{align*}
\int_0^t  \|g(s, \;\cdot\;)\|_{\Gamma, N}\, ds &= \int_0^t  \frac{(1+s)^{1+\varepsilon}}{(1+s)^{1+\varepsilon}} \|g(s, \;\cdot\;)\|_{\Gamma, N}\, ds \\
&\leqq \sup_{s \in [0, t]} (1+s)^{1+\varepsilon} \|g(s, \;\cdot\;)\|_{\Gamma, N} \int_0^t  \frac{ds}{(1+s)^{1+\varepsilon}}
\end{align*}
And since $$C:=\sup_{t \geq 0}\int_0^t  \frac{ds}{(1+s)^{1+\epsilon}} < \infty,$$ taking the supremum over time $t \geqq 0$ on both sides of the inequality yields the desired result.
\end{proof}

\begin{proposition}\label{prop:3} Let $u$ be a solution of ~\eqref{eq:LKG} satisfying ~\eqref{eq:IVP'} and ~\eqref{2.3}. Then, with constant $C>0$ and $\varepsilon$ a small positive number, we have
\begin{align*}
    |u(t, x)| \leqq C(1+t)^{-5 / 4}\left[\Delta_{5}+E_{1+\varepsilon, 5}(g)\right]
\end{align*}
for every $t \geqq 0, x \in \mathbb{R}^{3}$ and $\Delta_{5}=\left.\left(\|u(t)\|_{\Gamma, 5}\right)\right|_{t=0}$, depending only on the initial data $u_{0}, u_{1}$.
\end{proposition}
\begin{proof}
We have, using Proposition ~\ref{prop:3-1} with $N=5$, 
\begin{align*}
    \|u\|_{5} \leqq C_1\left(E_{1+\varepsilon, 5}(g)+\Delta_{5}\right)
\end{align*}
Also note that, from Proposition ~\ref{prop:1}, we have
\begin{align*}
    |u(t, x)| \leqq C_2(1+t)^{-5 / 4}\left(\|u\|_{5}+E_{1+\varepsilon, 3}(g)\right)
\end{align*}
We obtain
\begin{align*}
    |u(t, x)| &\leqq C_2(1+t)^{-5 / 4}\left(C_1\left(E_{1+\varepsilon, 5}(g)+\Delta_{5}\right)+E_{1+\varepsilon, 3}(g)\right) \\
    &\leqq C(1+t)^{-5 / 4}(\Delta_{5}+E_{1+\varepsilon, 5}(g))
\end{align*}
\end{proof}

\section{Energy Estimates for the Perturbed Klein-Gordon Equation}\label{section:hormander}

Consider the perturbed Klein-Gordon equation
\begin{equation}\label{PKG}\tag{P.K.G.}
\square_1 u + \sum_{j, k=0}^3 \gamma^{jk} \cdot \partial_j \partial_k u = f,
\end{equation}
where $\partial_0 = -\partial_t$, and $\{\gamma^{jk}:\mathbb{R}_+\times\mathbb{R}^3\rightarrow\mathbb{R} \mid j,k\in \{0..3\}\}$ are smooth functions of space-time. We subject ~\eqref{PKG} to the initial conditions ~\eqref{eq:IVP'}, and aim to use the energy method to bound $\mathcal{E}(u(t, \;\cdot\;))$  using L. Hörmander's method outlined in \cite{hormander1997}. This will be done in theorem ~\eqref{hormander}. To do so, we first prove a few useful results:
\begin{lemma}[Gronwall's Inequality]\label{lemma:gronwall}
Let $u \in C(\mathbb{R}_+, \mathbb{R})$ satisfying the differential inequality
\[
u'(t) \leqq a(t) u(t) + b(t) \quad \text{for all} \quad t\in\mathbb{R}_+,
\]
for some $a, b \in L^1(\mathbb{R}_+)$. Then, $u$ satisfies the pointwise estimate
\[
u(t) \leqq e^{A(t)} u(0) + \int_0^t b(s) e^{A(t)-A(s)} \, ds, \quad \text{for all} \quad t \in \mathbb{R}_+,
\]
where we have defined the primitive function
\[
A(t) := \int_0^t a(s) \, ds.
\]
Some examples and important special cases of the Gronwall lemma are
\begin{align*}
& u'(t) \leqq a(t) u(t) \quad \Longrightarrow \quad u(t) \leqq u(0) e^{A(t)}, \\
& u'(t) \leqq a u(t) + b \quad \Longrightarrow \quad u(t) \leqq u(0) e^{at} + \frac{b}{a} \left(e^{at} - 1\right), \\
& u'(t) \leqq a u(t) + b(t) \quad \Longrightarrow \quad u(t) \leqq u(0) e^{at} + \int_0^t e^{a(t-s)} b(s) \, ds, \\
& u'(t) + b(t) \leqq a(t) u(t), \quad\text{with}\quad a, b \geqq 0 \quad \Longrightarrow \quad u(t) + \int_0^t b(s) \, ds \leqq u(0) e^{A(t)}.
\end{align*}
\end{lemma}

\begin{proof}
Define $\phi(t) := e^{-A(t)} u(t)$ and compute its derivative:
\[
\phi'(t) = -a(t)e^{-A(t)}u(t) + e^{-A(t)}u'(t).
\]
Substitute $u'(t) \leq a(t)u(t) + b(t)$ into $\phi'(t)$:
\[
\phi'(t) \leq -a(t)e^{-A(t)}u(t) + e^{-A(t)}(a(t)u(t) + b(t)) = e^{-A(t)}b(t).
\]
Integrate both sides from $0$ to $t$:
\[
\int_0^t \phi'(s) \, ds \leq \int_0^t e^{-A(s)}b(s) \, ds.
\]
The LHS becomes $\phi(t) - \phi(0)$, thus
\[
\phi(t) \leq \phi(0) + \int_0^t e^{-A(s)}b(s) \, ds.
\]
Recall $\phi(t) = e^{-A(t)}u(t)$ and $\phi(0) = u(0)$:
\[
e^{-A(t)}u(t) \leq u(0) + \int_0^t e^{-A(s)}b(s) \, ds.
\]
Multiply both sides by $e^{A(t)}$ to obtain the pointwise estimate:
\[
u(t) \leq e^{A(t)}u(0) + \int_0^t b(s) e^{A(t)-A(s)} \, ds.
\]
\end{proof}

\begin{lemma}\label{boxudell0u}
For any smooth function $u:\mathbb{R}_+\times\mathbb{R}^3\rightarrow\mathbb{R}$, 
\begin{align*}
\square u \cdot \partial_{0} u = \frac{1}{2} \partial_{0} \left| \partial_{0} u \right|^2 + \frac{1}{2} \partial_{0} \sum_{i=1}^{3} \left| \partial_{i} u \right|^{2} - \sum_{i=1}^{3} \partial_{i} \left( \partial_{0} u \cdot \partial_{i} u \right)
\end{align*}
where $\partial_0 = -\partial_t$.
\end{lemma}
\begin{proof}
\begin{align*}
\square u \cdot \partial_{0} u &= \partial_{0}^2 u \cdot \partial_{0} u - \sum_{i=1}^{3} \partial_{i}^2 u \cdot \partial_{0} u \\
&= \partial_{0}^2 u \cdot \partial_{0} u + \sum_{i=1}^{3} \partial_{i} u \cdot \partial_{0} \partial_{i}u - \sum_{i=1}^{3} \partial_{i} u \cdot \partial_{0} \partial_{i}u - \sum_{i=1}^{3} \partial_{i}^2 u \cdot \partial_{0} u \\
&= \frac{1}{2} \partial_{0} \left| \partial_{0} u \right|^2 + \frac{1}{2} \partial_{0} \sum_{i=1}^{3} \left| \partial_{i} u \right|^2 - \sum_{i=1}^{3} \partial_{i} \left( \partial_{0} u \cdot \partial_{i} u \right).
\end{align*}
\end{proof}

\begin{lemma}\label{thm:div-1}
For any smooth function $u:\mathbb{R}_+\times\mathbb{R}^3\rightarrow\mathbb{R}$, and for any $j,k\in\{0..3\}$,
\begin{align*}
\int_{\mathbb{R}^3} \partial_{i} \left(\partial_{0} u \cdot \partial_{j} u\right) \, dx = 0
\end{align*}
where $\partial_0 = -\partial_t$.
\end{lemma}
\begin{proof}
Define the vector field \(V\) such that \(V_k = 0\) for \(k \neq i\) and \(V_i = \partial_{0} u \cdot \partial_{j} u\). Therefore, the divergence of \(V\), \(\text{div}(V) = \partial_{i} \left(\partial_{0} u \cdot \partial_{j} u\right)\), encapsulates the directed derivative of the product of \(\partial_{0} u\) and \(\partial_{j} u\) in the \(i\)-th direction. By applying the divergence theorem (lemma ~\eqref{thm:div}) and considering that the divergence of \(V\) over the entirety of \(\mathbb{R}^3\) equates to the flux through the boundary, which vanishes at infinity, the integral simplifies to \(0\), completing the proof.
\end{proof}

\begin{lemma}\label{dell0udelljdellku}
For any smooth function $u:\mathbb{R}_+\times\mathbb{R}^3\rightarrow\mathbb{R}$, and for any $j,k\in\{1..3\}$,
\begin{align*}
    \partial_0 u \cdot \partial_j \partial_k u = \frac{1}{2}\partial_j(\partial_0 u \cdot \partial_k u) + \frac{1}{2}\partial_k(\partial_0 u \cdot \partial_j u) - \frac{1}{2}\partial_0(\partial_j u \cdot \partial_k u).
\end{align*}
\end{lemma}
\begin{proof}
Applying Leibniz rule, compute:
\begin{align*}
    \partial_j(\partial_0 u \cdot \partial_k u) &= \partial_j \partial_0 u \cdot \partial_k u + \partial_0 u \cdot \partial_j \partial_k u \\
    \partial_k(\partial_0 u \cdot \partial_j u) &= \partial_k \partial_0 u \cdot \partial_j u + \partial_0 u \cdot \partial_k \partial_j u \\
    -\partial_0(\partial_j u \cdot \partial_k u) &= -\partial_0 \partial_j u \cdot \partial_k u - \partial_j u \cdot \partial_0 \partial_k u \\
\end{align*}
Taking their sum yields the desired result. 
\end{proof}

In the following, we introduce a notation for linear combinations which will be useful in proving energy methods, where we do not care about the constants that the operators $\Gamma$ are multiplied by, but we do care about the order of these operators. Moreover, from now on, multiplication by constants $C$ will often be used without regard to their different values. This means that their values can change from a line to another even though their symbol $C$ will not change. This is of course without any consequence on our results. 
\begin{lemma}\label{lemma:GammaAlphaPartial}
Given a 10-index $\alpha$ of size $|\alpha|=N$, the commutator of $\Gamma^{\alpha}$ with any first order partial derivative $\partial$ is a linear combination of terms of the form $\partial \Gamma^\beta$, with $\beta$ being a 10-index of size at most $N-1$. Using our linear combination notation, we write:
\begin{align*}
    [\Gamma^{\alpha}, \partial] = \lcomb{|\beta|\leqq N-1}{}{} \partial \Gamma^\beta
\end{align*}
\end{lemma}
\begin{proof}
    The proof is done by induction on the size $N$ of $\alpha$. For $N=1$, if $\Gamma^{\alpha} = \partial$, the commutator is immediately $0$ since partial derivatives commute. If $\Gamma^{\alpha} = \Omega_{ab}$, we have, for any $c\in\{0..3\}$,
    \begin{align*}
        [ \Omega_{ab}, \partial_c ] &= [ x_a \partial_b - x_b \partial_a, \partial_c ] = x_a [\partial_b, \partial_c] - x_b [\partial_a, \partial_c] + [\partial_b, \partial_c] x_a - [\partial_a, \partial_c] x_b \\
        &= x_a \delta_{bc} \partial - x_b \delta_{ac} \partial - x_a \delta_{bc} \partial + x_b \delta_{ac} \partial = 0
    \end{align*}
    where we used the Kronecker delta defined by:
    \begin{equation*}
    \delta_{ij} = 
    \begin{cases} 
    1 & \text{if } i = j,\\
    0 & \text{otherwise}.
    \end{cases}
    \end{equation*}
    Thus, for $N=1$, the commutator $[\Gamma^{\alpha}, \partial]$ is indeed a linear combination of terms $\partial \Gamma^\beta$ with $|\beta|\leqq N-1=0$, since the zero term can be considered as a linear combination of no $\Gamma^\beta$ terms. Now, assume the hypothesis to hold for some $N\in\mathbb{N}$ and let us show it holds for $N+1$. For $|\alpha| = N + 1$, write $\Gamma^{\alpha} = \Gamma^{i}\Gamma^{\beta}$, where $|\beta| = N$. Then,
    \begin{align*}
        [\Gamma^{\alpha}, \partial] &= [\Gamma^{i}\Gamma^{\beta}, \partial] = \Gamma^{i}[\Gamma^{\beta}, \partial] + [\Gamma^{i}, \partial]\Gamma^{\beta} \\
        &= \Gamma^{i} \left( \lcomb{|\gamma|\leqq N-1}{}{} \partial \Gamma^\gamma \right) + 0\cdot\Gamma^{\beta} \quad (\text{by induction hypothesis})\\
        &= \lcomb{|\gamma|\leqq N-1}{}{} \Gamma^{i} \partial \Gamma^\gamma = \lcomb{|\gamma|\leqq N-1}{}{} \partial \Gamma^{i} \Gamma^\gamma - [\Gamma^{i}, \partial]\Gamma^\gamma\\
        &= \lcomb{|\gamma|\leqq N}{}{} \partial \Gamma^\gamma
    \end{align*}
    where we combined $\Gamma^{i}$ and $\Gamma^\gamma$ into a $\Gamma$ operator of order 1 higher in the last line. This concludes the proof. 
\end{proof}

\begin{lemma}
Given a 10-index $\alpha$ of size $|\alpha|=N$, the commutator of $\Gamma^{\alpha}$ with any second order partial derivative $\partial \partial$ is a linear combination of terms of the form $\partial \partial \Gamma^\beta$, with $\beta$ being a 10-index of size at most $N-1$. Using our linear combination notation, we write:
\begin{align*}
    [\Gamma^{\alpha}, \partial \partial] = \lcomb{|\beta|\leqq N-1}{}{} \partial \partial \Gamma^\beta
\end{align*}
\end{lemma}
\begin{proof}
Using the previous lemma, compute 
    \begin{align*}
        [\Gamma^{\alpha}, \partial \partial] &= \Gamma^{\alpha} \partial \partial - \partial \partial \Gamma^{\alpha} = \Gamma^{\alpha} \partial \partial - \partial \Gamma^{\alpha} \partial - \partial [\Gamma^{\alpha}, \partial] = (\Gamma^{\alpha} \partial - \partial \Gamma^{\alpha}) \partial - \partial [\Gamma^{\alpha}, \partial] \\
        &= [\Gamma^{\alpha}, \partial] \partial - \partial [\Gamma^{\alpha}, \partial] = \left(\lcomb{|\beta|\leqq N-1}{}{} \partial \Gamma^\beta\right) \partial - \partial \left(\lcomb{|\beta|\leqq N-1}{}{} \partial \Gamma^\beta\right) \\
        &= \left(\lcomb{|\beta|\leqq N-1}{}{} \partial \Gamma^\beta \partial \right) + \left(\lcomb{|\beta|\leqq N-1}{}{} \partial \partial \Gamma^\beta\right) = \lcomb{|\beta|\leqq N-1}{}{} \partial \partial \Gamma^\beta - \lcomb{|\beta|\leqq N-1}{}{} \partial [\Gamma^\beta, \partial] \\
        &= \lcomb{|\beta|\leqq N-1}{}{} \partial \partial \Gamma^\beta + \lcomb{|\beta|\leqq N-2}{}{} \partial \partial \Gamma^\beta = \lcomb{|\beta|\leqq N-1}{}{} \partial \partial \Gamma^\beta
    \end{align*}
\end{proof}

\begin{theorem}[Hörmander's Lemma 7.4.1]\label{hormander}
Let \(u\) be a global solution of the perturbed Klein-Gordon equation ~\eqref{PKG}. If \(u\) vanishes for large \(|x|\), and if
\[
\sum_{j, k=0}^3\:\sup_{\substack{t\geq0 \\ x\in\mathbb{R}^3}}\:\left|\gamma^{jk}(t,x)\right| \leqq \frac{1}{2},
\]
then, for any \(t\geqq0\),
\[
\mathcal{E}(u(t,\;\cdot\;)) \leqq C\left(\mathcal{E}(u(0, \;\cdot\;)) + \int_0^t \|f(s, \;\cdot\;)\|_{L^2} \, ds\right) \exp \left(\int_0^t C \Gamma(s) \, ds\right),
\]
where for all $t\geqq0$, 
\[
\Gamma(t) = \sum_{i, j, k=0}^3\:\sup_{x\in\mathbb{R}^3}\:\left|\partial_i \gamma^{jk}(t, x)\right|,
\]
and $\mathcal{E}^2(u(t, \;\cdot\;))$ is defined as in ~\eqref{def:energynorm}. 
\end{theorem}

\begin{proof}
Define
\begin{align*} 
\mathcal{E}_+^2(u(t, \;\cdot\;)) = \mathcal{E}^2(u(t, \;\cdot\;)) + \int_{\mathbb{R}^3} \left[\gamma^{00} (\partial_0 u)^2 - \sum_{j, k=1}^3 \gamma^{jk} \partial_j u \, \partial_k u\right] \, dx,
\end{align*}
and note that for any $t\geqq0$, 
\begin{align*} 
\mathcal{E}(u(t, \;\cdot\;)) \leqq \mathcal{E}_+(u(t, \;\cdot\;)).
\end{align*}
Let us show that, for any $t\geqq0$,
\begin{align*} 
\mathcal{E}(u(t, \;\cdot\;)) \leqq C\left(\mathcal{E}(u(0, \;\cdot\;)) + \int_0^t \|f(s, \;\cdot\;)\|_{L^2} ds\right) \exp{\left(C\int_0^t\Gamma(s)ds\right)}
\end{align*}
To do so, we multiply both sides of \eqref{PKG} by \(\partial_0 u\) and then integrate over \(\mathbb{R}^3\) term by term. Note the use of Lemmas ~\eqref{thm:div-1} and ~\eqref{boxudell0u} in the following computations:
\begin{align*}
\int_{\mathbb{R}^3} \square_1 u \cdot \partial_0 u \, dx &= \int_{\mathbb{R}^3} \square u \cdot \partial_0 u \, dx + \int_{\mathbb{R}^3} u \cdot \partial_0 u \, dx \\
&= \frac{1}{2} \partial_{0} \int_{\mathbb{R}^3} \left|\partial_{0} u\right|^2 \, dx + \frac{1}{2} \partial_{0} \sum_{i=1}^{3} \int_{\mathbb{R}^3} \left|\partial_{i} u\right|^2 \, dx \\
& \quad \quad - \sum_{i=1}^{3} \int_{\mathbb{R}^3} \partial_{i} (\partial_{0} u \cdot \partial_{i} u) \, dx + \frac{1}{2} \partial_0 \int_{\mathbb{R}^3} |u|^2 \, dx\\
&= \frac{1}{2} \partial_0 \left\| \partial_0 u(t, \;\cdot\;) \right\|_{L^2}^2 + \frac{1}{2} \partial_0 \sum_{i=1}^{3} \left\| \partial_i u(t, \;\cdot\;) \right\|_{L^2}^2 + 0 + \frac{1}{2} \partial_0 \|u(t, \;\cdot\;)\|_{L^2}^2\\
&= \frac{1}{2} \partial_0 \mathcal{E}^2(u(t, \;\cdot\;))
\end{align*}
Now, using Lemma ~\eqref{dell0udelljdellku}, we aim to get rid of terms of the form $\partial \partial u$ since we cannot bound them by energy norms: 
\begin{align*}
\sum_{j, k=0}^3 \gamma^{jk} \cdot \partial_j \partial_k u \cdot \partial_0 u
&= \frac{1}{2} \sum_{j, k=0}^3 \gamma^{jk} \cdot \partial_j(\partial_0 u \cdot \partial_k u) + \frac{1}{2} \sum_{j, k=0}^3 \gamma^{jk} \cdot \partial_k(\partial_0 u \cdot \partial_j u) \\
& \quad - \frac{1}{2} \sum_{j, k=0}^3 \gamma^{jk} \cdot \partial_0(\partial_j u \cdot \partial_k u) \\
&= \frac{1}{2} \partial_0(\gamma^{00} \cdot |\partial_0 u|^2) + \frac{1}{2} \sum_{j, k=1}^3 \partial_j(\gamma^{jk} \cdot \partial_0 u \cdot \partial_k u) \\ 
& \quad - \frac{1}{2} \partial_0 \gamma^{00} \cdot |\partial_0 u|^2 - \frac{1}{2} \sum_{j, k=1}^3 \partial_j \gamma^{jk} \cdot \partial_0 u \cdot \partial_k u \\ 
& \quad + \frac{1}{2} \partial_0(\gamma^{00} \cdot |\partial_0 u|^2) + \frac{1}{2} \sum_{j, k=1}^3 \partial_k(\gamma^{jk} \cdot \partial_0 u \cdot \partial_j u) \\
& \quad - \frac{1}{2} \partial_0 \gamma^{00} \cdot |\partial_0 u|^2 - \frac{1}{2} \sum_{j, k=1}^3 \partial_k \gamma^{jk} \cdot \partial_0 u \cdot \partial_j u \\ 
& \quad - \frac{1}{2} \partial_0(\gamma^{00} \cdot |\partial_0 u|^2) - \frac{1}{2} \sum_{j, k=1}^3 \partial_0(\gamma^{jk} \cdot \partial_j u \cdot \partial_k u) \\ 
& \quad + \frac{1}{2} \partial_0 \gamma^{00} \cdot |\partial_0 u|^2 + \frac{1}{2} \sum_{j, k=1}^3 \partial_0 \gamma^{jk} \cdot \partial_k u \cdot \partial_j u 
\end{align*}
We end up with terms of the form $\partial u$. We can now integrate over $\mathbb{R}^3$, making use of Lemma ~\eqref{thm:div-1} to let some terms vanish:
\begin{align}
\int_{\mathbb{R}^3}\sum_{j, k=0}^3 \gamma^{jk} \cdot & \partial_j \partial_k u \cdot \partial_0 u \, dx
= \frac{1}{2} \int_{\mathbb{R}^3} \partial_0 \left[\gamma^{00} \cdot |\partial_0 u|^2 - \sum_{j, k=1}^3 \gamma^{jk} \cdot \partial_j u \cdot \partial_k u \right] \, dx \nonumber \\
& \quad - \frac{1}{2} \int_{\mathbb{R}^3} \partial_0 \gamma^{00} \cdot |\partial_0 u|^2 \, dx - \frac{1}{2} \sum_{j, k=1}^3 \int_{\mathbb{R}^3} \partial_j \gamma^{jk} \cdot \partial_0 u \cdot \partial_k u \, dx \nonumber \\
& \quad - \frac{1}{2} \sum_{j, k=1}^3 \int_{\mathbb{R}^3} \partial_k \gamma^{jk} \cdot \partial_0 u \cdot \partial_j u \, dx + \frac{1}{2} \sum_{j, k=1}^3 \int_{\mathbb{R}^3} \partial_0 \gamma^{jk} \cdot \partial_k u \cdot \partial_j u \, dx \nonumber \\
& \quad + \frac{1}{2} \sum_{j, k=1}^3 \int_{\mathbb{R}^3} \partial_j(\gamma^{jk} \cdot \partial_0 u \cdot \partial_k u) \, dx + \frac{1}{2} \sum_{j, k=1}^3 \int_{\mathbb{R}^3} \partial_k(\gamma^{jk} \cdot \partial_0 u \cdot \partial_j u) \, dx \nonumber \\
&= \frac{1}{2} \int_{\mathbb{R}^3} \partial_0 \left[\gamma^{00} \cdot |\partial_0 u|^2 - \sum_{j, k=1}^3 \gamma^{jk} \cdot \partial_j u \cdot \partial_k u \right] \, dx \label{eq:a} \\
& \quad - \frac{1}{2} \sum_{j, k=1}^3 \int_{\mathbb{R}^3} \partial_k \gamma^{jk} \cdot \partial_0 u \cdot \partial_j u \, dx + \frac{1}{2} \sum_{j, k=1}^3 \int_{\mathbb{R}^3} \partial_0 \gamma^{jk} \cdot \partial_k u \cdot \partial_j u \, dx \label{eq:b} \\
& \quad + \frac{1}{2} \sum_{j, k=1}^3 \int_{\mathbb{R}^3} \partial_j(\gamma^{jk} \cdot \partial_0 u \cdot \partial_k u) \, dx + \frac{1}{2} \sum_{j, k=1}^3 \int_{\mathbb{R}^3} \partial_k(\gamma^{jk} \cdot \partial_0 u \cdot \partial_j u) \, dx \label{eq:c}
\end{align}
We keep the terms of line ~\eqref{eq:a} in the LHS and move the terms of lines ~\eqref{eq:b} and ~\eqref{eq:c} to the RHS. Note that all terms in lines ~\eqref{eq:b} and ~\eqref{eq:c} can be simplified to a term of the form $\lcomb{j,k}{}{} \int_{\mathbb{R}^3} \partial \gamma^{jk} \cdot \partial u \cdot \partial u \, dx$. We obtain:
\begin{align*}
\frac{1}{2}\partial_0 \left[\mathcal{E}_+^2(u(t, \;\cdot\;))\right] &= \int_{\mathbb{R}^3} f \cdot \partial_0 u \, dx + \lcomb{j,k}{}{} \int_{\mathbb{R}^3} \partial \gamma^{jk} \cdot \partial u \cdot \partial u \, dx \\
&\leqq C\left(\|f(t, \;\cdot\;)\|_{L^2} \cdot \|\partial_0 u \|_{L^2} + \Gamma(t) \cdot \|\partial u \|_{L^2} \cdot \|\partial u \|_{L^2} \right)\\
&\leqq C\left(\|f(t, \;\cdot\;)\|_{L^2} \cdot \mathcal{E}(u(t, \;\cdot\;)) + \Gamma(t) \cdot \mathcal{E}^2(u(t, \;\cdot\;)) \right)\\
&\leqq C\left(\|f(t, \;\cdot\;)\|_{L^2} \cdot \mathcal{E}_+(u(t, \;\cdot\;)) + \Gamma(t) \cdot \mathcal{E}_+^2(u(t, \;\cdot\;)) \right)
\end{align*}
i.e.,
\begin{align*}
\partial_0 \mathcal{E}_+(u(t, \;\cdot\;)) &\leqq C\left(\|f(t, \;\cdot\;)\|_{L^2} + \Gamma(t) \cdot \mathcal{E}_+(u(t, \;\cdot\;)) \right)
\end{align*}
And we make use of Gronwall's lemma ~\eqref{lemma:gronwall} to get:
\begin{align*}
\mathcal{E}_+(u(t, \;\cdot\;)) &\leqq C\left( \mathcal{E}_+(u(0, \;\cdot\;)) + \int_0^t\|f(s, \;\cdot\;)\|_{L^2}ds \right) \exp{\left(C\int_0^t \Gamma(s)ds\right)} 
\end{align*}
Now, assuming the following result to be true for all $t\geqq0$:
\begin{equation}\label{key}
|\mathcal{E}_+^2(u(t, \;\cdot\;)) - \mathcal{E}^2(u(t, \;\cdot\;))| \leqq \frac{1}{2} \mathcal{E}^2(u(t, \;\cdot\;)),
\end{equation}
We get:
\begin{equation*}
\mathcal{E}^2(u(t, \;\cdot\;))\leqq \frac{3}{2} \mathcal{E}_0^2(u(t, \;\cdot\;))
\end{equation*}
\begin{equation*}
\mathcal{E}_+(u(t, \;\cdot\;))\leqq C \mathcal{E}(u(t, \;\cdot\;))
\end{equation*}
And since $\mathcal{E}(u(t, \;\cdot\;))\leqq \mathcal{E}_+(u(t, \;\cdot\;))$, we get the desired result. 
So it remains to show ~\eqref{key}:
\begin{align*}
|\mathcal{E}_+^2(u(t, \;\cdot\;)) &- \mathcal{E}^2(u(t, \;\cdot\;))| \leqq \int_{\mathbb{R}^3} \left( \left|\gamma^{00} \cdot (\partial_0 u)^2 \right| + \sum_{j, k=1}^3 \left| \gamma^{jk} \cdot \partial_j u \cdot \partial_k u\right| \right) \, dx \\
&\leqq \sum_{j, k=0}^3 \int_{\mathbb{R}^3} \left| \gamma^{jk} \cdot \partial_j u \cdot \partial_k u\right| \, dx \leqq \sum_{j, k=0}^3\:\sup_{\substack{t\geq0 \\ x\in\mathbb{R}^3}}\:\left|\gamma^{jk}(t,x)\right| \int_{\mathbb{R}^3}\left| \partial_j u \cdot \partial_k u\right| \, dx \\
&\leqq \frac{1}{2} \sum_{j, k=0}^3 \|\partial_j u \|_{L^2} \cdot \|\partial_k u \|_{L^2} \leqq \frac{1}{2}\mathcal{E}^2(u(t, \;\cdot\;))
\end{align*}
\end{proof}

\section{Proof of the Theorem}\label{section:proof}
For simplicity, and without much less generality, we assume that the nonlinear term $F$ in \eqref{eq:NKG} is affine in $u^{\prime\prime}$, i.e.,
\begin{equation*}
F\left(u, u^{\prime}, u^{\prime \prime}\right)=\sum_{a, i = 0}^3 f^{a i}\left(u, u^{\prime}\right) \partial_{a i}^2 u+\sum_{a=0}^3 f^a\left(u, u^{\prime}\right) \partial_a u
\end{equation*}
where $f^{a i}, f^a$ are smooth functions of $u, u^{\prime}$ with $f^{a i}(0,0)=f^a(0,0)=0$ for all $a=0,1,2,3, i=1,2,3$. 
We may assume that
\begin{equation*}
\sum_{a, i = 0}^3\left|f^{a i}\left(u, u^{\prime}\right)\right| \leqq \frac{1}{2}
\end{equation*}
for any $u, u^{\prime}$ with $|u|+\left|u^{\prime}\right| \leqq 1$. For the sake of simplicity, we will often use $f^j$ to refer to either $f^a$ or $f^{ai}$. Thus, for any $j$, $f^j \in \{f^{a}, f^{ai}\mid a,i\in \{0..3\}\}$. 

The proof of \eqref{thm:1} follows a standard iteration argument. Consider 
a smooth function $w_0 : \mathbb{R}_+ \times \mathbb{R}^3 \rightarrow \mathbb{R}$ which satisfies the following conditions:
\begin{subequations}
\begin{align}
\textit{Finite speed of propagation:}& \quad w_0(t, x)=0 \quad \text{for any } t \geq 0, \quad|x| \geq t+1, \label{cond:a}\\ 
\textit{Decay at infinity:}& \quad |w_0|_{5 / 4,8} \leq K \varepsilon, \label{cond:b}\\ 
\textit{Control over the energy:}& \quad \|w_0\|_{14} \leq K \varepsilon, \label{cond:c}
\end{align}
\end{subequations}
for some fixed $K$, sufficiently large, and $\varepsilon$ an arbitrary small positive constant.

We construct a sequence of functions $w_n : \mathbb{R}_+\times\mathbb{R}^3\rightarrow\mathbb{R}$ defined recursively as follows: 
for any $n \in \mathbb{N}$, $w_{n+1}$ is the solution to the linear inhomogeneous Klein-Gordon equation
\begin{equation}\label{eq:4}
\Box_1 w_{n+1} = \sum_{a, i = 0}^3 f^{ai}\left(w_n, w_n^{\prime}\right) \cdot \partial_a \partial_i w_{n+1} +\sum_{a=0}^3 f^a\left(w_n, w_n^{\prime}\right) \cdot \partial_a w_{n+1},
\end{equation}
satisfying the initial data \eqref{eq:IVP'}. The existence of solutions to such equations is a classical result (see \cite{polyanin2002handbook}). Assuming that for all $n \in \mathbb{N}$, $w_n$ satisfies  conditions \eqref{cond:a}, \eqref{cond:b} and \eqref{cond:c}, we can easily show by induction that $(w_n)$ is a Cauchy sequence in the Hilbert space $H^{14}$ and thus converges to a limit $w_\infty \in H^{14}$, since $H^{14}$ is complete. Thus, it is clear that $w_\infty$ is a global solution of the Klein-Gordon equation:
\begin{equation}
\Box_1 w_\infty = \sum_{a, i = 0}^3 f^{ai}\left(w_\infty, w_\infty^{\prime}\right) \cdot \partial_a \partial_i w_\infty +\sum_{a=0}^3 f^a\left(w_\infty, w_\infty^{\prime}\right) \cdot \partial_a w_\infty,
\end{equation}
which proves the main theorem of this paper. Thus it remains to show that for all $n \in \mathbb{N}$, $w_n$ satisfies  conditions \eqref{cond:a}, \eqref{cond:b} and \eqref{cond:c}. The proof is done by induction. The base case is fulfilled at $n=0$ by the assumptions on $w_0$. It remains to prove the inductive step. For the sake of simplicity, we fix $n \in \mathbb{N}$ and set $u:=w_n$, $v:=w_{n+1}$, and we show the following:

\begin{theorem}\label{thm:1'}
Let $u$ be a smooth function satisfying \eqref{cond:a}, \eqref{cond:b} and \eqref{cond:c}, and let $v$ be the unique solution of the problem:
\begin{equation}\label{thm1-eq}
\Box_1 v = \sum_{a, i = 0}^3 f^{ai}\left(u, u^{\prime}\right) \cdot \partial_a \partial_i v +\sum_{a=0}^3 f^a\left(u, u^{\prime}\right) \cdot \partial_a v,
\end{equation}
satisfying \eqref{eq:IVP'}. Let $\Delta_N = \Delta_N(u_0, u_1)$ be defined as in \eqref{def:Delta}. We claim that if $K$ is sufficiently large, independent of $\varepsilon$, and $\varepsilon$ is sufficiently small, and
\begin{equation}\label{cond:d}
    \Delta_{14}(u_0, u_1) \leqq \varepsilon,
\end{equation}
then $v$ satisfies the same conditions \eqref{cond:a}, \eqref{cond:b}, and \eqref{cond:c} as $u$.
\end{theorem}

To simplify notation we shall denote from now on by $|\,\cdot\,|_{N}$ the norm $|\,\cdot\,|_{5 / 4, N}$ introduced in \eqref{def:normkN}. 
To prove ~\eqref{thm:1'}, we make use of the following results:

\begin{lemma}\label{lemma:GammaBetafj}
    For any $f^j \in \{f^{a}, f^{ai}\mid a,i\in \{0..3\}\}$ and for any multi-index $\beta$ of size $|\beta| = N$,
    \begin{equation*}
        \left\| \Gamma^{\beta} f^{j}(u, u^\prime) \right\|_{L^2} \leqq C \left\| u \right\|_{\Gamma, N}
    \end{equation*}
\end{lemma}
\begin{proof}
    Recalling that $f^j(0,0) = 0$, begin by writing the Taylor series of $f^j$ centered at $0$:
    $$
    f^j(u, u^\prime) = \sum_{n\geqq1} \: \sum_{k=0}^n \: C_{nk}^j \: u^k \cdot (u^\prime)^{n-k}
    $$
    Then apply $\Gamma^\beta$ to get:
    $$
    \Gamma^\beta f^j(u, u^\prime) = \sum_{n\geqq1} \: \sum_{k=0}^n \: \lcomb{|\beta_1| + |\beta_2| \leqq N} \: \Gamma^{\beta_1} [u^k] \cdot \Gamma^{\beta_2} [(u^\prime)^{n-k}]
    $$
    Then, taking the $L^2$-norm on both sides, we make use of lemma ~\eqref{lemma:GammaAlphaPower}, and change $C$ at every line, we get:
    \begin{align*}
        \left\|\Gamma^\beta f^j(u, u^\prime)\right\|_{L^2} &\leqq \sum_{n\geqq1} \: \sum_{k=0}^n \: \lcomb{|\beta_1| + |\beta_2| \leqq N} \left\| \Gamma^{\beta_1} [u^k] \right\|_{L^2} \cdot \left\| \Gamma^{\beta_2} [(u^\prime)^{n-k}] \right\|_{L^2} \\
        &\leqq \sum_{n\geqq1} \: \sum_{k=0}^n \: \lcomb{|\gamma_1| + |\gamma_2| \leqq N} \left\| \Gamma^{\gamma_1} u \right\|_{L^2}^k \cdot \left\| \Gamma^{\gamma_2} u^\prime \right\|_{L^2}^{n-k}\\
        &\leqq C \sum_{n\geqq1} \: \sum_{k=0}^n \: \left\| u \right\|_{\Gamma, N}^k \cdot \left\| u \right\|_{\Gamma, N}^{n-k}\\
        &\leqq C \left\| u \right\|_{\Gamma, N} \sum_{n\geqq0} \: \sum_{k=0}^n \: \left\| u \right\|_{\Gamma, N}^n\\
        &\leqq C \left\| u \right\|_{\Gamma, N}
    \end{align*}
\end{proof}

\begin{lemma}
    For any $f^j \in \{f^{a}, f^{ai}\mid a,i\in \{0..3\}\}$ and for any multi-index $\beta$ of size $|\beta| = N$,
    \begin{equation*}
        \left| \Gamma^{\beta} f^{j}(u, u^\prime) \right|_{L^\infty} \leqq C \left| u \right|_{\Gamma, N}
    \end{equation*}
\end{lemma}
\begin{proof}
    The proof is the same as that of lemma ~\eqref{lemma:GammaBetafj}, replacing all the $L^2$-norms by $L^\infty$-norms. 
\end{proof}

\begin{proposition}[Generalized energy estimates]\label{3.5}
    Let $v$ be the solution of ~\eqref{thm1-eq} satisfying ~\eqref{eq:IVP'}. For any integer $N>0$, there exists a constant $C_{N}>0$ such that
    \begin{equation}\label{eq:3.5}
        \|v\|_{N+1} \leqq C_{N}\left(\Delta_{N+1}+|v|_{\lceil(N+1) / 2\rceil} \cdot \|u\|_{N+1}\right) \cdot \exp \left\{C_{N}|u|_{\lceil(N+1) / 2\rceil}\right\}
    \end{equation}
\end{proposition}
\begin{proof}
The proof of follows precisely the pattern of proving energy estimates for perturbed Klein-Gordon equations, as in the proof for ~\eqref{hormander}. Let $N\in \mathbb{N}$, let $\alpha$ be a multi-index such that $|\alpha|\leqq N+1$, and let us show that:  
\begin{equation}\label{EGammaAlpha}
    \partial_0\mathcal{E}(\Gamma^\alpha v(t, \;\cdot\;)) \leqq C \left(|u(t, \;\cdot\;)|_{\Gamma, \left\lceil\frac{N+1}{2}\right\rceil} \cdot \mathcal{E}(\Gamma^\alpha v(t, \;\cdot\;)) + \|u(t, \;\cdot\;)\|_{\Gamma, N+1} \cdot |v(t, \;\cdot\;)|_{\Gamma, \left\lceil\frac{N+1}{2}\right\rceil}\right)
\end{equation}
To do so, we start by applying the operator $\Gamma^\alpha$ to equation ~\eqref{thm1-eq}.  By making use of the commutation properties of the operator $\Gamma^\alpha$  with $\Box_{1}$ (see lemma ~\eqref{lemma:GammaAlphaBox1}), we obtain:
\begin{equation*}
\Box_1 \Gamma^\alpha v -  \Gamma^\alpha \sum_{a, i = 0}^3 f^{ai}\left(u, u^{\prime}\right) \cdot \partial_a \partial_i v = \Gamma^\alpha \sum_{a=0}^3 f^a\left(u, u^{\prime}\right) \cdot \partial_a v
\end{equation*}
Now, compute:
\begin{align*}
\Gamma^\alpha \sum_{a, i = 0}^3 & f^{ai}\left(u, u^{\prime}\right) \cdot \partial_a \partial_i v = \sum_{a, i = 0}^3 f^{ai}\left(u, u^{\prime}\right) \cdot \Gamma^\alpha \partial_a \partial_i v \\
& \quad \quad + \sum_{a,i=0}^3 \:\lcomb{|\alpha_1|+|\alpha_2|\leqq N} \Gamma^{\alpha_1} f^{ai}(u, u^\prime) \cdot  \Gamma^{\alpha_2} \partial_a \partial_i v \\
&= \sum_{a, i = 0}^3 f^{ai}\left(u, u^{\prime}\right) \cdot \partial_a \partial_i \Gamma^\alpha v + \sum_{a,i=0}^3 \:\lcomb{|\beta|\leqq N} f^{ai}(u, u^\prime) \cdot \partial \partial \Gamma^{\beta} v \\
&\quad + \sum_{a,i=0}^3 \:\lcomb{|\alpha_1|+|\alpha_2|\leqq N} \Gamma^{\alpha_1} f^{ai}(u, u^\prime) \cdot  \partial_a \partial_i \Gamma^{\alpha_2} v \\
& \quad \quad + \sum_{a,i=0}^3 \:\lcomb{|\alpha_1|+|\alpha_2|\leqq N-1} \Gamma^{\alpha_1} f^{ai}(u, u^\prime) \cdot  \partial \partial \Gamma^{\alpha_2} v \\
&= \sum_{a, i = 0}^3 f^{ai}\left(u, u^{\prime}\right) \cdot \partial_a \partial_i \Gamma^\alpha v + \sum_{a,i=0}^3 \:\lcomb{|\beta|\leqq N+1} f^{ai}(u, u^\prime) \cdot \partial \Gamma^{\beta} v \\
&\quad + \sum_{a,i=0}^3 \:\lcomb{|\alpha_1|+|\alpha_2|\leqq N+1} \Gamma^{\alpha_1} f^{ai}(u, u^\prime) \cdot  \partial \Gamma^{\alpha_2} v \\
& \quad \quad + \sum_{a,i=0}^3 \:\lcomb{|\alpha_1|+|\alpha_2|\leqq N} \Gamma^{\alpha_1} f^{ai}(u, u^\prime) \cdot  \partial \Gamma^{\alpha_2} v \\
&= \sum_{a,i=0}^3 f^{ai}(u, u^\prime) \cdot \partial \partial \Gamma^\alpha v + \sum_{a,i=0}^3 \:\lcomb{|\alpha_1|+|\alpha_2|\leqq N+1} \Gamma^{\alpha_1} f^{ai}(u, u^\prime) \cdot \partial \Gamma^{\alpha_2} v
\end{align*}
We keep the terms $\sum_{a,i=0}^3 f^{ai}(u, u^\prime) \cdot \partial \partial \Gamma^\alpha v$ in the LHS and move the other terms to the RHS. Similarly, compute:
\begin{align*}
\Gamma^\alpha & \sum_{a = 0}^3 f^{a}\left(u, u^{\prime}\right) \cdot \partial_{a} v = \sum_{a = 0}^3 \:\lcomb{|\alpha_1|+|\alpha_2| \leqq N+1} \Gamma^{\alpha_1}f^{a}\left(u, u^{\prime}\right) \cdot \Gamma^{\alpha_2}\partial_{a} v \\
&= \sum_{a = 0}^3 \:\lcomb{|\alpha_1|+|\beta| \leqq N+1} \Gamma^{\alpha_1}f^{a}\left(u, u^{\prime}\right) \cdot \partial_{a} \Gamma^{\alpha_2} v + \sum_{a = 0}^3 \:\lcomb{|\alpha_1|+|\alpha_2| \leqq N} \Gamma^{\alpha_1}f^{a}\left(u, u^{\prime}\right) \cdot \partial \Gamma^{\beta} v \\
&= \sum_{a=0}^3 \:\lcomb{|\alpha_1|+|\alpha_2|\leqq N+1} \Gamma^{\alpha_1} f^{a}(u, u^\prime) \cdot \partial \Gamma^{\alpha_2} v
\end{align*}
Re-indexing the functions $\{f^{a}, f^{ai}\mid a,i\in \{0..3\}\}$ with an index $j$, and putting it all together, equation ~\eqref{thm1-eq} becomes:
\begin{equation*}
\Box_1 \Gamma^\alpha v - \sum_{j} f^{j}(u, u^\prime) \cdot \partial \partial \Gamma^\alpha v = \sum_{j} \:\lcomb{|\alpha_1|+|\alpha_2|\leqq N+1} \Gamma^{\alpha_1} f^{j}(u, u^\prime) \cdot \partial \Gamma^{\alpha_2} v := f
\end{equation*}
We now aim to bound the $L^2$-norm of the RHS in order to follow the same steps as in the proof for Hörmander's lemma (theorem ~\eqref{hormander}):
\begin{align*}
    \left\|f(t, \;\cdot\;)\right\|_{L^2} &\leqq \sum_{j} \lcomb{|\alpha_1|+|\alpha_2|\leqq N+1} \left\| \Gamma^{\alpha_1} f^{j}(u, u^\prime) \cdot \partial \Gamma^{\alpha_2} v\right\|_{L^2} \\
    &\leqq \sum_{j} \lcomb{\substack{|\beta_1|\leqq \left\lceil\frac{N+1}{2}\right\rceil \\ |\beta_2| \leqq N+1}} \left| \Gamma^{\beta_1} f^{j}(u, u^\prime) \right|_{L^\infty} \cdot \left\|\partial \Gamma^{\beta_2} v\right\|_{L^2} \\
    & \quad \quad + \sum_{j} \lcomb{\substack{|\beta_1|\leqq N+1 \\ |\beta_2| \leqq \left\lceil\frac{N+1}{2}\right\rceil}} \left\| \Gamma^{\beta_1} f^{j}(u, u^\prime) \right\|_{L^2} \cdot \left|\partial \Gamma^{\beta_2} v\right|_{L^\infty}\\
    &\leqq C\left( \left| 
    u(t, \;\cdot\;) \right|_{\Gamma, \left\lceil\frac{N+1}{2}\right\rceil} \cdot \left\| v(t, \;\cdot\;)\right\|_{\Gamma, N+1} + \left\| u(t, \;\cdot\;) \right\|_{\Gamma, N+1} \cdot \left|v(t, \;\cdot\;)\right|_{\Gamma, \left\lceil\frac{N+1}{2}\right\rceil} \right) 
\end{align*}
Thus, using the exact same steps as in the proof for Hörmander's lemma, we get ~\eqref{EGammaAlpha}. 
Now, we can sum both sides of inequality ~\eqref{EGammaAlpha} over all multi-indices $\alpha$ size $|\alpha| \leqq N+1$ to get:
\begin{align*}
    \partial_0 &\|v(t, \;\cdot\;)\|_{\Gamma, N+1} =  \sum_{|\alpha| \leq N+1} \partial_0\mathcal{E}(\Gamma^\alpha v(t, \;\cdot\;)) \\
    &\leq C \left(|u(t, \;\cdot\;)|_{\Gamma, \left\lceil\frac{N+1}{2}\right\rceil} \sum_{|\alpha| \leq N+1} \mathcal{E}(\Gamma^\alpha v(t, \;\cdot\;)) + \|u(t, \;\cdot\;)\|_{\Gamma, N+1} \cdot |v(t, \;\cdot\;)|_{\Gamma, \left\lceil\frac{N+1}{2}\right\rceil}\right) \\
    &\leq C \left(|u(t, \;\cdot\;)|_{\Gamma, \left\lceil\frac{N+1}{2}\right\rceil} \cdot \|v(t, \;\cdot\;)\|_{\Gamma, N+1} + \|u(t, \;\cdot\;)\|_{\Gamma, N+1} \cdot |v(t, \;\cdot\;)|_{\Gamma, \left\lceil\frac{N+1}{2}\right\rceil}\right)
\end{align*}
Now, applying Gronwall's inequality (lemma ~\eqref{lemma:gronwall}) yields the desired result.
\end{proof}

\begin{proposition}[Decay estimates]\label{prop:3.6}
Let $v$ be the solution of ~\eqref{thm1-eq} satisfying ~\eqref{eq:IVP'}. For every $M \in \mathbb{N}$, there exists a constant $C_M\geq0$ such that
\begin{equation}\label{eq:3.6}
    |v|_{M} \leqq C_{M}\left(\Delta_{M+5}+|u|_{\lceil(M+6) / 2\rceil} \cdot\|v\|_{M+6}+|v|_{\lceil(M+6) / 2\rceil} \cdot\|u\|_{M+6}\right)
\end{equation}
\end{proposition}

\begin{proof}
Let $M\in\mathbb{N}$ and let $\alpha$ be a multi-index of size $|\alpha|\leqq M$. We apply $\Gamma^\alpha$ to both sides of ~\eqref{thm1-eq} to get:
\begin{align*}
\Box_1 \Gamma^\alpha v &= \Gamma^\alpha \sum_{a, i = 0}^3 f^{ai}\left(u, u^{\prime}\right) \cdot \partial_a \partial_i v + \Gamma^\alpha \sum_{a=0}^3 f^a\left(u, u^{\prime}\right) \cdot \partial_a v \\
&= \sum_{a, i = 0}^3 \;\lcomb{|\alpha_1|+|\alpha_2|\leqq M+1}\; \Gamma^{\alpha_1} f^{ai}\left(u, u^{\prime}\right) \cdot \Gamma^{\alpha_2} \partial v \\
& \quad + \sum_{a, i = 0}^3 \;\lcomb{|\alpha_1|+|\alpha_2|\leqq M}\; \Gamma^{\alpha_1} f^{a}\left(u, u^{\prime}\right) \cdot \Gamma^{\alpha_2} \partial_a v \\
&= \sum_{j} \;\lcomb{|\alpha_1|+|\alpha_2|\leqq M+1}\; \Gamma^{\alpha_1} f^{j}\left(u, u^{\prime}\right) \cdot \Gamma^{\alpha_2} \partial v := g \\
\end{align*}
Calling the RHS $g$, we have, using Proposition ~\eqref{prop:3}, 
\begin{align*}
        |\Gamma^\alpha v(t, x)| \leqq C(1+t)^{-5 / 4}\left[\Delta_{M+5}+E_{1+\varepsilon, 5}(g)\right]. 
\end{align*}
So let us compute $E_{1+\varepsilon, 5}(g)$. Given a multi-index $\beta$ of size $|\beta|\leqq 5$, we have:
\begin{align*}
    \left\|\Gamma^\beta g(t, \;\cdot\;)\right\|_{L^2} &\leqq \sum_{j} \;\lcomb{|\alpha_1|+|\alpha_2|\leqq M+6}\; \left\|\Gamma^{\alpha_1} f^{j}\left(u, u^{\prime}\right) \cdot \Gamma^{\alpha_2} \partial v \right\|_{L^2} \\
    &\leqq \sum_{j} \;\lcomb{\substack{|\beta_1|\leqq \left\lceil\frac{M+6}{2}\right\rceil \\ |\beta_2|\leqq M+6}}\; \left|\Gamma^{\beta_1} f^{j}\left(u, u^{\prime}\right) \right|_{L^\infty} \cdot \left\|\Gamma^{\beta_2} \partial v \right\|_{L^2} \\
    & \quad + \sum_{j} \;\lcomb{\substack{|\beta_1|\leqq M+6 \\ |\beta_2|\leqq \left\lceil\frac{M+6}{2}\right\rceil}}\; \left\|\Gamma^{\beta_1} f^{j}\left(u, u^{\prime}\right) \right\|_{L^2} \cdot \left|\Gamma^{\beta_2} \partial v \right|_{L^\infty} \\
    &\leqq C \left(|u(t, \;\cdot\;)|_{\Gamma, \left\lceil\frac{M+6}{2}\right\rceil} \cdot \|v(t, \;\cdot\;)\|_{\Gamma, M+6} + \|u(t, \;\cdot\;)\|_{\Gamma, M+6} \cdot |v(t, \;\cdot\;)|_{\Gamma, \left\lceil\frac{M+6}{2}\right\rceil} \right) 
\end{align*}
Thus, 
\begin{align*}
    E_{1+\varepsilon, 5}(g) = \sup_{t\geqq0} (1+t)^{1+\varepsilon} \|g(t, \;\cdot\;)\|_{\Gamma, 5} = \sup_{t\geqq0} (1+t)^{1+\varepsilon} \sum_{|\beta|\leqq4} \|g(t, \;\cdot\;)\|_{L^2} \\
    \leqq C \left(|u|_{\left\lceil\frac{M+6}{2}\right\rceil} \cdot \|v\|_{M+6} + \|u\|_{M+6} \cdot |v|_{\left\lceil\frac{M+6}{2}\right\rceil} \right)
\end{align*}
\end{proof}

\begin{proof}[Proof of Theorem ~\eqref{thm:1'}]
The vanishing of $v$ for $|x| \geqq t+1$ follows immediately from the properties of $u, u_{0}, u_{1}$ if $f^{a i} \equiv 0$, and is in general a consequence of sharp uniqueness theorems of the type proved by F. John in \cite{fjohn}.
Now, taking $N=13$ in ~\eqref{eq:3.5}, we have:
\begin{align*}
\|v\|_{14} \leq C_{13} \left(\Delta_{14} + |v|_{7} \cdot \|u\|_{14}\right) \cdot \exp \left(C_{13}|u|_7\right)
\end{align*}
and using ~\eqref{cond:b}, ~\eqref{cond:c} together with ~\eqref{cond:d}, we infer that
\begin{align*}
\|v\|_{14} \leq C_{13}\left(\varepsilon+\varepsilon K|v|_{7}\right) \exp \left(\varepsilon K C_{14}\right)
\end{align*}
Picking $\varepsilon>0$ sufficiently small such that
\begin{equation}\label{epsilon-cond-1}
    \varepsilon K C_{14} \leqq 1,
\end{equation}
we derive
\begin{equation}\label{3.7'}
    \|v\|_{14} \leqq 3 \varepsilon C_{13} \left(1 + K|v|_{7}\right)
\end{equation}
On the other hand, using ~\eqref{eq:3.6} with $M=8$ we see that
\begin{align*}
    |v|_{8} \leqq C_{8}\left(\Delta_{13}+|u|_{7} \cdot\|v\|_{14}+|v|_{7} \cdot\|u\|_{14}\right)
\end{align*}
And again using ~\eqref{cond:b}, ~\eqref{cond:c} together with ~\eqref{cond:d}, get
\begin{align*}   
    |v|_{8} &\leqq C_{8}\left(\varepsilon+K\varepsilon \|v\|_{14}+|v|_{7} K\varepsilon\right) \\
    &\leqq C_{8}\varepsilon\left(1+K\varepsilon 3 C_{13} \left(1 + K|v|_{7}\right)+|v|_{7} K\right) \\
    &\leqq C_{8}\varepsilon\left(1+|v|_{7} K(1 + K\varepsilon 3 C_{13}) + K\varepsilon 3 C_{13}\right) \\
    &\leqq C_{8}\varepsilon\left(4+4|v|_{7} K\right) \\
    &\leqq 4 C_{8}\varepsilon\left(1+|v|_{8} K\right)
\end{align*}
Hence, if
\begin{align}\label{epsilon-cond-2}
    4 C_8 \varepsilon K \leqq \frac{1}{2}
\end{align}
we obtain
$$
|v|_{8} \leqq\left(8 C_{8}\right) \varepsilon
$$
which proves ~\eqref{cond:b} provided that
\begin{align}\label{K-cond-1}
    K \geqq 8 C_{8}
\end{align}
Finally, from ~\eqref{3.7'},
\begin{align*}  
    \|v\|_{14} &\leqq 3 \varepsilon C_{13} \left(1 + K|v|_{7}\right) \\
    &\leqq 3 \varepsilon C_{13} \left(1 + K^2\varepsilon\right) \\
    &\leqq K\varepsilon 
\end{align*}
provided that
\begin{align}\label{K-cond-2}
    K \geqq 6C_{13},
\end{align}
and,
\begin{align}\label{epsilon-cond-3}
    \varepsilon K^{2} \leqq 1,
\end{align}
which proves ~\eqref{cond:c}. 
Thus, picking $K$ by ~\eqref{K-cond-1}, ~\eqref{K-cond-2} and the picking $\varepsilon$ sufficiently small by ~\eqref{epsilon-cond-1}, ~\eqref{epsilon-cond-2}, and ~\eqref{epsilon-cond-3}, we have proved Theorem ~\eqref{thm:1'}.
\end{proof}

\begin{appendices}

\section{$L^p$ spaces and Sobolev spaces}
In this section we introduce basic objects and results fundamental to this paper: multi-indices, Sobolev spaces and norms, Sobolev inequalities. See \cite{evans2010} for more details. 
\begin{definition}[\(L^p\) norm]
Let \(0<p<\infty\) and let \((X, \mathcal{M}, \mu)\) denote a measure space. If \(f: X \rightarrow \mathbb{R}\) is a measurable function, then we define
\[
\|f\|_{L^p(X)}:=\left(\int_X|f|^p \, d\mu\right)^{\frac{1}{p}} \quad \text{and} \quad |f|_{L^{\infty}(X)}:=\operatorname{ess\,sup}_{x \in X}|f(x)|.
\]
Note that \(\|f\|_{L^p(X)}\) may take the value \(\infty\).
\end{definition}

\begin{definition}[The space \(L^p(X)\)]
The space \(L^p(X)\) is the set
\[
L^p(X)=\left\{f: X \rightarrow \mathbb{R} \mid \|f\|_{L^p(X)} < \infty\right\}.
\]
The space \(L^p(X)\) satisfies the following vector space properties:
\begin{enumerate}
    \item For each \(\alpha \in \mathbb{R}\), if \(f \in L^p(X)\) then \(\alpha f \in L^p(X)\);
    \item If \(f, g \in L^p(X)\), then
    \[
    |f+g|^p \leq 2^{p-1}\left(|f|^p+|g|^p\right),
    \]
    so that \(f+g \in L^p(X)\).
    \item The triangle inequality is valid if \(p \geq 1\).
\end{enumerate}
The most interesting cases are \(p=1,2, \infty\), while all of the \(L^p\) spaces arise often in nonlinear estimates.
\end{definition}

\begin{theorem} \(L^{p}(\Omega)\) is a Banach space
\end{theorem}

\begin{definition}[Test functions]
For $\Omega \subset \mathbb{R}^n$, set
\[
C_0^{\infty}(\Omega)=\left\{u \in C^{\infty}(\Omega) \mid u \;\text{has compact support}\right\},
\]
the smooth functions with compact support.
\end{definition}

\begin{definition}[Multi-index]\label{def:multi-index}
An element $\alpha \in \mathbb{Z}^n$ is called an $n$-index. For such an $\alpha= \left(\alpha_1, \ldots, \alpha_n\right)$, we write $D^\alpha=\frac{\partial^{\alpha_1}}{\partial x_{\alpha_1}} \cdots \frac{\partial^{\alpha_n}}{\partial x_{\alpha_n}}$ and we define the size of the multi-index by $|\alpha|=\alpha_1+\cdots+ \alpha_n$.
\end{definition}

\begin{definition}[Weak derivative]
Suppose that $u \in L_{\text{loc}}^1(\Omega)$. Then $v^\alpha \in L_{\text{loc}}^1(\Omega)$ is called the $\alpha^{\text{th}}$ weak derivative of $u$, written $v^\alpha=D^\alpha u$, if for any $\phi \in C_0^{\infty}(\Omega)$,
\[
\int_{\Omega} u(x) D^\alpha \phi(x) \, dx=(-1)^{|\alpha|} \int_{\Omega} v^\alpha(x) \phi(x) \, dx.
\]
\end{definition}

\begin{remark}
Note that if the weak derivative exists, it is unique. To see this, suppose that both $v_1$ and $v_2$ are the weak derivative of $u$ on $\Omega$. Then $\int_{\Omega}(v_1-v_2) \phi \, dx=0$ for all $\phi \in C_0^{\infty}(\Omega)$, so that $v_1=v_2$ a.e.
\end{remark}

\begin{definition}[Sobolev space \(W^{k, p}(\Omega)\)]
For integers $k \geq 0$ and $1 \leq p \leq \infty$,
\[
W^{k, p}(\Omega)=\left\{u \in L_{\mathrm{loc}}^1(\Omega) \mid D^\alpha u \text{ exists and is in } L^p(\Omega) \text{ for }|\alpha| \leq k\right\}.
\]
\end{definition}

\begin{definition}[Norm in \(W^{k, p}(\Omega)\)]\label{def:sobolev-norms}
For $u \in W^{k, p}(\Omega)$ define
\[
\|u\|_{W^{k, p}(\Omega)}=\left(\sum_{|\alpha| \leq k}\left\|D^\alpha u\right\|_{L^p(\Omega)}^p\right)^{\frac{1}{p}} \text{ for } 1 \leq p<\infty
\]
and
\[
|u|_{W^{k, \infty}(\Omega)}=\max_{|\alpha| \leq k}\left|D^\alpha u\right|_{L^{\infty}(\Omega)}.
\]
It is clear that these functions define norms since they are a finite sum of $L^p$ norms.
\end{definition}

\begin{theorem}\label{lemma:complete} \(W^{k, p}(\Omega)\) is a Banach space
\end{theorem}
\begin{proof}
Let \(u_j\) denote a Cauchy sequence in \(W^{k, p}(\Omega)\). It follows that for all \(|\alpha| \leq k\), \(D^\alpha u_j\) is a Cauchy sequence in \(L^p(\Omega)\). Since \(L^p(\Omega)\) is a Banach space, for each \(\alpha\) there exists \(u^\alpha \in L^p(\Omega)\) such that
\[
D^\alpha u_j \rightarrow u^\alpha \text{ in } L^p(\Omega).
\]
When \(\alpha=(0, \ldots, 0)\) we set \(u:=u^{(0, \ldots, 0)}\) so that \(u_j \rightarrow u\) in \(L^p(\Omega)\). We must show that \(u^\alpha=D^\alpha u\).
For each \(\phi \in C_0^{\infty}(\Omega)\),
\[
\begin{aligned}
\int_{\Omega} u D^\alpha \phi \, dx & = \lim_{j \rightarrow \infty} \int_{\Omega} u_j D^\alpha \phi \, dx \\
& = (-1)^{|\alpha|} \lim_{j \rightarrow \infty} \int_{\Omega} D^\alpha u_j \phi \, dx \\
& = (-1)^{|\alpha|} \int_{\Omega} u^\alpha \phi \, dx,
\end{aligned}
\]
thus, \(u^\alpha=D^\alpha u\) and hence \(D^\alpha u_j \rightarrow D^\alpha u\) in \(L^p(\Omega)\) for each \(|\alpha| \leq k\), which shows that \(u_j \rightarrow u\) in \(W^{k, p}(\Omega)\).
\end{proof}

\begin{definition}[Hilbert space \(H^k(\Omega)\)]
For integers \(k \geq 0\) and \(p=2\), we define
\[
H^k(\Omega) = W^{k, 2}(\Omega).
\]
\(H^k(\Omega)\) is a Hilbert space with inner-product \((u, v)_{H^k(\Omega)} = \sum_{|\alpha| \leq k}(D^\alpha u, D^\alpha v)_{L^2(\Omega)}\).
\end{definition}

\begin{lemma}\label{lemma:injection}
Let $\Omega$ be a bounded open subset of $\mathbb{R}^n$ with $C^1$ boundary, and let $u\in H^k(\Omega)$. If $k>n/p$, then
\begin{equation*}
|u|_{L^\infty(\Omega)} \leqq C\|u\|_{H^k \left(\Omega\right)}
\end{equation*}
with the constant $C$ depending only on $n, k, \Omega$
\end{lemma}

\section{Additional Useful Results }
In this section, we show some additional results that are fundamental for the proof of this paper. Most of these results are not directly used in the main proof and are very computational, hence the choice of grouping them in a separate appendix. 
\begin{lemma}[Lorentz invariance of $\Box$]\label{lemma:lorentz-invatiance}
    For any $\Gamma^i \in \Gamma$, 
    \begin{equation*}
        [\Gamma^i, \Box] = 0
    \end{equation*}
\end{lemma}
\begin{proof}
    When $\Gamma^i = \partial$, the result is straightforward since partial derivatives commute. When $\Gamma^i = \Omega_{ab}$, we compute:
    \begin{align*}
        [\Omega_{ab}, \Box] = [\Omega_{ab}, \partial_t^2 - \Delta] = [\Omega_{ab}, \partial_t^2] - [\Omega_{ab}, \Delta]
    \end{align*}
    And, since
    \begin{align*}
        [\Omega_{ab}, \partial^2] &= \Omega_{ab}\partial^2 - \partial^2\Omega_{ab} \\
        &= (x_a\partial_b - x_b\partial_a)\partial^2 - \partial^2(x_a\partial_b - x_b\partial_a) \\
        &= x_a\partial_b\partial^2 - x_b\partial_a\partial^2 - x_a\partial^2\partial_b + x_b\partial^2\partial_a \\
        &= 0,
    \end{align*} 
    we get the desired result. 
\end{proof}
\begin{lemma}
    For any $10$-index $\alpha$, for any fixed ordering of the operators $\Gamma$,
    \begin{equation*}
        [\Gamma^\alpha, \Box] = 0
    \end{equation*}
\end{lemma}   
\begin{proof}
    The proof is done by (strong) induction on the size $N:=|\alpha|$ of the multi-index. The base case ($N=1$) is given by the previous lemma. Assume the lemma holds for any $10$-index $\beta$ of size $|\beta|\leqq N$ and let $\alpha$ be a $10$-index of size $|\alpha|=N+1$. Then we can write $\Gamma^\alpha = \Gamma^i \Gamma^\beta$ where $\beta$ is a $10$-index of size $|\beta|=N$ and $\Gamma^i \in \Gamma$. Then we compute: 
    \begin{align*}
        [\Gamma^\alpha, \Box] &= [\Gamma^i \Gamma^\beta, \Box] = \Gamma^i \Gamma^\beta \Box - \Box \Gamma^i \Gamma^\beta = \Gamma^i \Box \Gamma^\beta - \Gamma^i \Box \Gamma^\beta = 0
    \end{align*}
\end{proof}
\begin{corollary}\label{lemma:GammaAlphaBox1}
    For any $10$-index $\alpha$, for any fixed ordering of the operators $\Gamma$,  
    \begin{equation*}
        [\Gamma^\alpha, \Box_1] = 0
    \end{equation*}
\end{corollary}
\begin{proof}
    This is a straightforward consequence of the previous lemma, the fact the fact that $[\Gamma^\alpha, 1] = 0$, and the bilinearity of the Lie bracket. 
\end{proof}

\begin{lemma}\label{lemma:commutator-of-Gammas}
    The commutator of any two $\Gamma$ operators is a $\mathbb{R}$-linear combination of the $\Gamma$ operators, i.e., for any $\Gamma^i,\Gamma^j \in \Gamma$, 
    \begin{align*}
        [\Gamma^i,\Gamma^j] \in \operatorname{Span_{\mathbb{R}}}(\Gamma)
    \end{align*}
\end{lemma}
\begin{proof}
When $\Gamma^i$ and $\Gamma^j$ are both partial derivative operators, their commutator is $0$ so the lemma holds. Consider the case $\Gamma^i = \Omega_{ab} = x_a \partial_b - x_b \partial_a$ and $\Gamma^j = \partial_c$. Then we have
\begin{align*}
    [\Omega_{ab}, \partial_c] &= [x_a \partial_b - x_b \partial_a, \partial_c] \\
    &= x_a [\partial_b, \partial_c] - x_b [\partial_a, \partial_c] + [\partial_c, x_a] \partial_b - [\partial_c, x_b] \partial_a \\
    &= 0 - 0 + (\delta_{ca} \partial_b - \delta_{cb} \partial_a) \\
    &= \delta_{ca} \partial_b - \delta_{cb} \partial_a,
\end{align*}
which is a $\mathbb{R}$-linear combination of elements of $\Gamma$. Now consider the case $\Gamma^i = \Omega_{ab} = x_a \partial_b - x_b \partial_a$ and $\Gamma^j = \Omega_{cd} = x_c \partial_d - x_d \partial_c$. Then we have:
\begin{align*}
    [\Omega_{ab}, \Omega_{cd}] &= [x_a \partial_b - x_b \partial_a, x_c \partial_d - x_d \partial_c] \\
    &= x_a [\partial_b, x_c \partial_d] - x_a [\partial_b, x_d \partial_c] - x_b [\partial_a, x_c \partial_d] + x_b [\partial_a, x_d \partial_c] \\
    &= x_a (\delta_{bc} \partial_d + x_c [\partial_b, \partial_d]) - x_a (\delta_{bd} \partial_c + x_d [\partial_b, \partial_c]) \\
    &\quad - x_b (\delta_{ac} \partial_d + x_c [\partial_a, \partial_d]) + x_b (\delta_{ad} \partial_c + x_d [\partial_a, \partial_c]) \\
    &= x_a \delta_{bc} \partial_d - x_a \delta_{bd} \partial_c - x_b \delta_{ac} \partial_d + x_b \delta_{ad} \partial_c,
\end{align*}
which shows that $[\Omega_{ab}, \Omega_{cd}]$ is also a $\mathbb{R}$-linear combination of $\Gamma$ operators and ends the proof. 
\end{proof}

\begin{lemma}\label{lemma:div-prod}
    Given a vector $\mathbf{F} : \mathbb{R}^3 \rightarrow \mathbb{R}^3$ and a scalar field $g : \mathbb{R}^3 \rightarrow \mathbb{R}$, we have
    \begin{align*}
        \operatorname{div}(g \cdot \mathbf{F}) = \nabla g \cdot \mathbf{F} + g \cdot \operatorname{div}(\mathbf{F})
    \end{align*}
    where \(\operatorname{div} \, \mathbf{F} = \nabla \cdot \mathbf{F}\) is the divergence of \(\mathbf{F}\). 
\end{lemma}
\begin{proof}
Let $\mathbf{F} = (F_1, F_2, F_3)$ be the components of the vector field and $g$ the scalar field. The divergence of $g\mathbf{F}$ is computed as
\begin{align*}
    \operatorname{div}(g\mathbf{F}) &= \nabla \cdot (g\mathbf{F}) \\
    &= \frac{\partial}{\partial x}(gF_1) + \frac{\partial}{\partial y}(gF_2) + \frac{\partial}{\partial z}(gF_3) \\
    &= \left(\frac{\partial g}{\partial x}F_1 + g\frac{\partial F_1}{\partial x}\right) + \left(\frac{\partial g}{\partial y}F_2 + g\frac{\partial F_2}{\partial y}\right) + \left(\frac{\partial g}{\partial z}F_3 + g\frac{\partial F_3}{\partial z}\right) \\
    &= (\nabla g \cdot \mathbf{F}) + g(\nabla \cdot \mathbf{F}),
\end{align*}
\end{proof}

\begin{lemma}[Divergence Theorem]\label{thm:div}
Let \(\Omega \subseteq \mathbb{R}^3\) be a region with smooth boundary \(\partial \Omega\), and let \(\mathbf{F}\) be a continuously differentiable vector field on \(\overline{\Omega}\). Then,
\begin{equation*}
\int_{\Omega} \operatorname{div} \, \mathbf{F} \, dV = \int_{\partial \Omega} \mathbf{F} \cdot \mathbf{n} \, dS,
\end{equation*}
where \(\operatorname{div} \, \mathbf{F} = \nabla \cdot \mathbf{F}\) is the divergence of \(\mathbf{F}\), \(\mathbf{n}\) is the outward unit normal vector on \(\partial \Omega\), \(dV\) is the volume element, and \(dS\) is the surface element.
\end{lemma}

\begin{lemma}\label{lemma:energy-helper}
    For any sufficiently smooth $u: \mathbb{R}_+ \times \mathbb{R}^3 \rightarrow \mathbb{R}$, 
    \begin{align*}
        \int_{\mathbb{R}^3}\partial_t u \cdot \Box_1 u \, dx = \frac{1}{2}\partial_t\mathcal{E}^2(u(t, \;\cdot\;))
    \end{align*}
    \end{lemma}
\begin{proof}
    \begin{align*}
        \int_{\mathbb{R}^3}\partial_t u \cdot \Box_1 u \, dx &= \int_{\mathbb{R}^3}\partial_t u \cdot \partial_t^2 u \, dx - \int_{\mathbb{R}^3}\partial_t u \cdot \Delta u \, dx + \int_{\mathbb{R}^3}\partial_t u \cdot u \, dx 
    \end{align*}
    We compute:
    \begin{align*}
        \int_{\mathbb{R}^3}\partial_t u \cdot \partial_t^2 u \, dx &= \int_{\mathbb{R}^3} \frac{1}{2}\partial_t \left(|\partial_t u|^2\right) \, dx = \frac{1}{2}\partial_t \left(\|\partial_t u(t, \;\cdot\;)\|_{L^2}^2\right) \\
        \int_{\mathbb{R}^3}\partial_t u \cdot u \, dx &= \int_{\mathbb{R}^3} \frac{1}{2}\partial_t \left(|u|^2\right) \, dx = \frac{1}{2}\partial_t \left(\| u(t, \;\cdot\;)\|_{L^2}^2\right) \\
    \end{align*}
    And, using lemmas ~\eqref{lemma:div-prod} and ~\eqref{thm:div}, compute:
    \begin{align*}
        \int_{\mathbb{R}^3} \partial_t u \cdot \Delta u \, dx &= \int_{\mathbb{R}^3} \operatorname{div}(\partial_t u \cdot \nabla u) \, dx + \int_{\mathbb{R}^3} \operatorname{div}(\partial_t \nabla u \cdot \nabla u) \, dx \\
        &= 0 + \int_{\mathbb{R}^3} \frac{1}{2}\partial_t \left(|\nabla u|^2\right) \, dx = \frac{1}{2}\partial_t \left(\| \nabla u(t, \;\cdot\;)\|_{L^2}^2\right)
    \end{align*}
\end{proof}

\begin{lemma}[Product rule for $\Gamma$ operators]
    For any $\Gamma^i \in \Gamma$ and for any sufficiently smooth functions $f, g: \mathbb{R}_+\times\mathbb{R}^3 \rightarrow \mathbb{R}$, we have
    \begin{align*}
    \Gamma^i(f \cdot g) = \Gamma^i f \cdot g + f \cdot \Gamma^i g.
    \end{align*}
\end{lemma}
\begin{proof}
    When $\Gamma^i$ is a partial derivative, the result follows directly from the product rule for derivatives:
    
    When $\Gamma^i = \Omega_{ab} = x_a \partial_b - x_b \partial_a$, we compute:
    \begin{align*}
    \Omega_{ab}(f \cdot g) &= (x_a \partial_b - x_b \partial_a)(f \cdot g) \\
    &= x_a \partial_b(f \cdot g) - x_b \partial_a(f \cdot g) \\
    &= x_a \left((\partial_b f) \cdot g + f \cdot (\partial_b g)\right) - x_b \left((\partial_a f) \cdot g + f \cdot (\partial_a g)\right) \\
    &= (x_a \partial_b f) \cdot g + (x_a f) \cdot (\partial_b g) - (x_b \partial_a f) \cdot g - (x_b f) \cdot (\partial_a g) \\
    &= \left(x_a \partial_b f - x_b \partial_a f\right) \cdot g + f \cdot \left(x_a \partial_b g - x_b \partial_a g\right) \\
    &= \Omega_{ab} f \cdot g + f \cdot \Omega_{ab} g.
    \end{align*}
\end{proof}

\begin{lemma}[Leibniz rule for $\Gamma$ operators]
    For any $10$-index $\alpha$ of size $|\alpha|=N$, for any $\Gamma^i \in \Gamma$ and for any sufficiently smooth $f, g: \mathbb{R}_+ \times \mathbb{R}^3 \rightarrow \mathbb{R}$, $\Gamma^\alpha(f \cdot g)$ is a linear combination of terms of the form $\Gamma^{\alpha_1} f \cdot \Gamma^{\alpha_2} g$, with $|\alpha_1| + |\alpha_2| \leqq N$. Using the notation introduced in lemma ~\eqref{lemma:GammaAlphaPartial}, we write:
    \begin{align*}
    \Gamma^\alpha(f \cdot g) = \lcomb{|\alpha_1| + |\alpha_2| \leqq N} \Gamma^{\alpha_1} f \cdot \Gamma^{\alpha_2} g 
    \end{align*}
\end{lemma}
\begin{lemma}\label{lemma:GammaAlphaPower}
    For any 10-index $\alpha$ of size $|\alpha|=N$, for any sufficiently smooth function $f :\mathbb{R}_+ \times \mathbb{R}^3 \rightarrow \mathbb{R}$, there exist a $10$-index $\beta$ of size $|\beta|\leqq N$ and $C\geqq0$ such that
    \begin{align*}
    \|\Gamma^\alpha(f^n)\| \leqq C\|\Gamma^\beta f\|^n
    \end{align*}
    where $f^n$ denotes the $n$-th power of $f$, and $\|\cdot\|$ denotes either the $L^2$-norm or the $L^\infty$-norm (taken over the spacial domain $\mathbb{R}^3$). 
\end{lemma}
\begin{proof}
    The proof is done by induction of the exponent $n$ of $f$. The base case ($n=0$) is trivial. Now assume the result holds for $n\in \mathbb{N}$, i.e., that for any $10$-index $\alpha_1$, there exists a $10$-index $\beta_{\alpha_1}$ of size $|\beta_{\alpha_1}|\leqq |\alpha_1|$ and $C\geqq0$ such that
    \begin{align*}
    \|\Gamma^{\alpha_1}(f^n)\| \leqq C\|\Gamma^{\beta_{\alpha_1}} f\|^n
    \end{align*}
    Let us show that it also holds for $n+1$. Let $\alpha$ be a 10-index of size $|\alpha| = N$ and we have, using the previous lemma and the hypothesis,
    \begin{align*}
        \|\Gamma^\alpha(f^{n+1})\| = \|\Gamma^\alpha(f^n \cdot f)\| &\leqq \lcomb{|\alpha_1| + |\alpha_2| \leqq N} \|\Gamma^{\alpha_1} (f^n)\| \cdot \|\Gamma^{\alpha_2} f \| \\ 
        &\leqq \lcomb{|\alpha_1| + |\alpha_2| 
        \leqq N} C\|\Gamma^{\beta_{\alpha_1}} f\|^n \cdot \|\Gamma^{\alpha_2} f \|
    \end{align*}
    Setting $\beta := \operatorname{argmax}_\gamma \|\Gamma^\gamma f\|$, where the dummy variable $\gamma$ ranges over the set $\{\beta_{\alpha_1} \mid |\alpha_1| \leqq N\} \cup \{\alpha_2 \mid |\alpha_2| \leqq N\}$. Then, changing $C$, we obtain
    \begin{align*}
        \|\Gamma^\alpha(f^{n+1})\| \leqq C\|\Gamma^\beta f\|^{n+1} , 
    \end{align*}
    which concludes the proof. 
\end{proof}

\end{appendices}

\section*{Acknowledgements}
I would like to thank my advisor Professor Annalaura Stingo for her patient guidance and support. I am also grateful for my academic training at the Bachelor Program of Ecole Polytechnique that provided me with the strong mathematical foundations to complete this thesis. 


\printbibliography

\end{document}